\documentclass[12pt,leqno]{article}
\usepackage{amsmath,amstext, amsfonts, verbatim,remreset}
\usepackage{graphicx}

\newtheorem {theorem} {Theorem}
\newtheorem {lemma} {Lemma}

\newcommand {\RR} {\mathbb R}
\newcommand {\NN} {\mathbb N}

\newcommand {\SSs} {\mathbb S}

\newcommand {\R} {\mathbb R}

\font\mathcal=cmbsy10 

\begin{document}

\title{Exact asymptotics of the optimal $L_{p,\Omega}$-error of linear spline interpolation}


\author{Vladislav Babenko, Yuliya Babenko\thanks{supported by SHSU Enhancement Grant for Research}, Dmytro Skorokhodov}
\date{}
\maketitle

{{ \centerline{\bf Abstract} In this paper we provide the exact
asymptotics of the optimal weighted $L_p$-error, $0<p< \infty$, of
linear spline interpolation of $C^2$ functions with positive
Hessian. The full description of the behavior of the optimal error
leads to the algorithm for construction of an asymptotically optimal
sequence of triangulations. In addition, we compute the minimum of
the $L_p$-error of linear interpolation of the function $x^2+y^2$
over all triangles of unit area for all $0<p<\infty$. This provides
the exact constant in the asymptotics of the optimal error. }}


\section{Definitions, history, and main results.}



\subsection {Definitions.}

Let $\RR^2$ be the space of points in the plane endowed with the
usual Euclidian distance. The distance between points $A$
and $B$ in $\RR^2$ we shall denote by $|AB|$.

Let $D=[0,1]^2 \subset \RR^2$. We use this region for simplicity;
the approach presented in this paper can be applied to any bounded
connected region which is a finite union of triangles. By $C(D)$ we
shall denote the space of functions continuous on $D$.

Let $L_{p}(D)$, $0<p\le \infty$, be the space of measurable functions $f:D\to\RR$ for which the value
$$
  \|f\|_p = \|f\|_{L_p(D)} : = \left\{\begin{array}{ll}
                \left(\displaystyle\int\limits_{D} |f (\tau) | ^p d\tau\right)^{\frac 1p},   &{\rm if}\;\;\;    0< p < \infty, \\ [10pt]
                {\rm ess sup} \{|f (\tau) |:\tau \in D \},                                    &{\rm if}\;\;\;       p =\infty.
\end{array}\right.
$$
is finite.

{\bf Remark:} Note that in the case $1\le p\le\infty$ the functional $\|\cdot\|_p$ is the usual norm in the space $L_p(D)$. For $p\in(0,1)$ the above relation defines a seminorm satisfying
\begin{equation}\label{norm_inverse}
  \|f+g\|_p^p\le\|f\|_p^p+\|g\|_p^p
\end{equation}
for arbitrary functions $f,g\in L_p(D)$.

Given a positive continuous weight function $\Omega \in C(D)$, we
define the weighted $L_p$-norm $\|\cdot\|_{p, \Omega}$ as
$$
\|f\|_{p, \Omega}=\|f\Omega\|_p.
$$
In the case $\Omega\equiv 1$, instead of $L_{p,\Omega}$ we shall simply
write $L_p$.

A collection $\triangle_N=\triangle_N(D)=\{T_i\}^N_{i=1}$ of
$N$ triangles in the plane is called a {\it triangulation} of a set
$D$ provided that
\begin{enumerate}
\item any pair of triangles from $\triangle_N$ intersect at most at
a common vertex or along a common edge,

\item $D=\displaystyle \bigcup\limits_{i=1}^N{T_i}$.
\end{enumerate}


Let ${\cal P}_1$ be the set of bivariate linear polynomials
$p(x,y)=ax+by+c$. Given a triangulation $\triangle_N$, define the
class of linear splines on $\triangle_N$ to be
$$
{\cal S}_1(\triangle_N):=\{ f \in C(D): \forall i=1,...,N\;\; \exists p\in P_1\;{\rm such\;that}\;f|_{T_i}=p|_{T_i}\}.
$$

Let
$s(f,\triangle_N)$ denote the spline from ${\cal S}_1(\triangle_N)$ which interpolates the function $f \in
C(D)$ at the vertices of the triangulation $\triangle_N$.
Note that the linear spline $s(f,\triangle_N)$ is uniquely determined
by its values at the vertices of the triangulation $\triangle_N$.

Now let the function $f\in C^2(D)$ and the number of triangles $N\in
\NN$ be fixed.
Set
$$   
R_N(f,L_{p,\Omega}):=\displaystyle \inf_{\triangle_N
}\|f-s(f,\triangle_N)\|_{p,\Omega},
$$   
where inf is taken over all triangulations of $D$ with $N$ triangles. This quantity will be called {\it optimal $L_{p,\Omega}$-error} of the interpolation
of the function $f$ by splines $s(f,\triangle_N)$.

 A triangulation $\triangle_N^0$ is called {\it optimal} for the given function $f$  if
$$   
\|f-s(f,\triangle^0_N)\|_{p,\Omega}=R_N(f,L_{p,\Omega}).
$$   

Exact values of $R_N(f,L_{p,\Omega})$ as well as optimal
triangulations $\triangle_N^0$ for every particular function $f$ can
be found only in few exceptional situations.
That is why the following two problems are interesting and important:

1) Find exact asymptotics of the optimal error $R_N(f,L_{p,\Omega})$ as $N \to \infty$ for any function  $f\in C^2(D)$.

2) Find an {\it asymptotically optimal} sequence of triangulations, i.e. a sequence of triangulations $\{ \triangle_N^* \}^{\infty}_{N=1}$ of $D$ such that

$$   
\displaystyle \lim_{N \to \infty}\frac{\|f-s(f,\triangle^*_N)\|_{p, \Omega} }{R_N(f,L_{p,\Omega})}=1.
$$   

For $f\in C^2(D)$ set $f_{xx}:=\frac{\partial^2f}{\partial x^2}$, $f_{xy}:=\frac{\partial^2f}{\partial x\partial y}$,
$f_{yy}:=\frac{\partial^2f}{\partial y^2}$, and denote the {\it Hessian} of $f$ by
$$
  H(f;x,y):=(f_{xx}f_{yy}-f_{xy}^2)(x,y).
$$

The main purpose of this paper is to solve the two formulated above
problems for functions $f\in C^2(D)$ such that $H(f;x,y)>0$ on $D$.
We shall need to resolve the following auxiliary extremal problem,
which is of independent interest.

For a given triangle $T\subset D$ and a function $g\in C(D)$ we
shall denote by $l(g,T)$ the linear polynomial which interpolates
$g$ at the vertices of the triangle $T$. In addition, for $0<
p<\infty$ and continuous weight function $\Omega$ set
$$
  d(g,T,L_{p,\Omega}):=\|g-l(g,T)\|_{L_{p,\Omega}(T)}.
$$
For $\overline{Q}(x,y):=x^2+y^2$ define
\begin{equation}\label{Cp}
  C_p^{+}:=\displaystyle \inf_{T} \frac{d(\overline{Q},T,L_p)}{ |T| ^{1+\frac 1p}},
\end{equation}
where $|T|$ denotes the area of the triangle $T$.

Main results of this paper are in the following two theorems.
\begin{theorem}\label{main_th}
Let $f \in C^2(D)$; $H(f;x,y)\ge C^+>0$ for all $(x,y) \in D$.
Let a positive continuous weight function
$\Omega(x,y)$ also be given. Then for all $0< p < \infty$
$$
\lim_{N \to \infty} N \cdot R_N(f,L_{p,\Omega})
  = \frac {C^+_p}{2} \|\sqrt{H}\Omega\|_{L_{\frac{p}{p+1}}(D)}  
$$
\end{theorem}

\begin{theorem}\label{T2}
For any $p>0$ infimum in the definition of $C_p^+$ is achieved only on
equilateral triangles. Consequently,
$$
  C_p^+=\left(\frac{4}{3\sqrt{3}}\right)^{1+\frac 1p}\int\limits_{T_0}(1-x^2-y^2)^p\,dx\,dy,
$$
where $T_0$ is the equilateral triangle with the center of
circumscribed circle at the origin.
\end{theorem}

{\bf Remark. }It is easy to see that
$$
  C_p^{+} = \displaystyle\left(\frac{4}{3\sqrt{3}}\right)^{1+\frac 1p}\left[\frac{\pi}{p+1}-6\int\limits_{1/2}^1 x(1-x^2)^p\arccos{\frac 1{2x}}\,dx\right]^{\frac 1p}.
$$
Moreover, if $B(a,b)$ is the Euler Beta function, and
$$
  B(x;a,b):=\int\limits_0^xt^{a-1}(1-t)^{b-1}\,dt
$$
is the incomplete Beta function, then
$$
  C_p^+=\displaystyle\left(\frac{4}{3\sqrt{3}}\right)^{1+\frac{1}{p}}\left[\frac{\pi}{p+1}-3{\rm B}\left(p+1,\frac{1}{2}\right){\rm B}\left(\frac{3}{4};p+\frac{3}{2},\frac{1}{2}\right)\right]^{\frac{1}{p}}.
$$

{\bf Corollary. \it Let $Q(x,y)=Ax^2+By^2+2Cxy$. Then for every $p\in(0,\infty)$
\begin{equation}\label{Cp1}
  \inf_{T}\frac{d(Q,T,L_p)}{|T|^{1+\frac{1}{p}}}=C_p^+\sqrt{AB-C^2}.
\end{equation}}

The idea of constructing an asymptotically optimal sequence of
triangulations is to substitute the function $f$ by piecewise
quadratic function $f_N$ for every $N$, for which the good
triangulation is constructed with the help of triangles solving
problem~(\ref{Cp1}).

In the proof of the lower estimate in Theorem~1 the fact that for an
arbitrary $f\in C^2(D)$ with positive Hessian there exists a
constant $\kappa$, depending on the function $f$ only, such that
$$
  d(f,T,L_p)\ge\frac{\kappa\cdot{\rm diam}\,T}{2^{5p+1}}\cdot|T|^{1+\frac
  1p},
$$
played an important role (see Lemmas~7,~8, and~10).

\subsection{History.}
The first result related to piecewise linear interpolation in two
dimensional case was obtained by L. Fejes Toth. He indicated
(\cite{Toth}, Ch. 5, \S 12) that for a body $C \subset \RR^3$ with
boundary of differentiability class $C^2$ and positive Gaussian
curvature $K(x,y)$ the Hausdorff distance of $C$ to its best
inscribed polytope with at most $n$ vertices is
$$
\frac{1+o(1)}{3\sqrt 3} \left ( \displaystyle \int _{\partial C} K(x,y)^{1/2} d\sigma(x,y) \right)\frac1n
$$
as $n\to \infty$, where $\sigma$ is the surface area measure on $\partial C$.
He also indicated that
the distance of $C$ to its best inscribed polytope with at most $n$ vertices (measured as a volume of the difference between $C$ and the polytope) is
$$
\frac{1+o(1)}{4\sqrt 3} \left ( \displaystyle \int _{\partial C} K(x,y)^{1/4} d\sigma(x,y) \right)^2\frac1n
$$
as $n\to \infty$. Even though all the ideas were mostly contained
in~\cite{Toth}, formally the complete proof was given by Gruber
(see~\cite{Gr2}). In addition, Gruber generalized these results to
higher dimensions, however the constants were implicit. He also
proved similar estimates for the error measured in symmetric
difference metric, Banach-Mazur metric, as well as using the
Schneider distance (see~\cite{Bor}). Among other interesting results
on these and closely related questions are results by
B$\rm{\ddot{o}}$r$\rm{\ddot{o}}$czky and Ludwig~\cite{Bor, KL}.
 Survey of further results on approximation of convex bodies by various polytopes in different metrics (inscribed, circumscribed, of the best approximation, with
  restrictions on the number of faces, etc.) can be found, for example, in~\cite{Bor, Gr2}.

With regard to the asymptotically optimal approximation of functions
by linear splines in different metrics the following results are
known.

Nadler in~\cite{Nadler} studied the sequence of asymptotically optimal triangulations
for the (in general discontinuous) piecewise linear approximation for an arbitrary $C^3$
function in the sense of minimizing error in the $L_2$-norm.

In the paper~\cite{us} this problem was solved for $p=\infty$ and $\Omega\equiv1$, and in~\cite{PhD}
it was extended for $p=\infty$ and positive weight functions $\Omega\in C(D)$.

Note that the case $p=\infty$, $\Omega\equiv1$ is close to (but not
indentical with) the result of L.~Fejes~Toth on approximation of
convex bodies by polygons in Hausdorff metric. The case $p=1$
follows from the results of B$\rm{\ddot{o}}$r$\rm{\ddot{o}}$czky and
Ludwig (see~\cite{KL}).


As for the question of computing the constant $C_p^+$ and the
optimality of the equilateral triangle, the following cases have
been investigated:
\begin{itemize}

\item $p=\infty$ (D'Azevedo and Simpson, in~\cite{Daz3}),

\item $p=1$ (B$\rm{\ddot{o}}$r$\rm{\ddot{o}}$czky and Ludwig, in~\cite{KL}),

\item $p=2$ (Pottmann, Hamann {\it et al}, in~\cite{kodla1}),

\item $p=3,4,5$ (Yu. Babenko, in~\cite{PhD}),

\item $p\in\NN$ (Chen, in~\cite{Chen1}).

\end{itemize}

The rest of the paper is organized as follows. Section 2 provides
certain preliminary results, in particular, on how affine
transformations affect the error of interpolation of a quadratic
function by linear splines. Section 3 contains the computation of
the constant $C_p^+$ in the case $1\le p<\infty$. The case
$p\in(0,1)$ is rather technical, and
thus is presented in the Appendix. 
In Section 4 we provide the proof of the upper estimate in
Theorem~1. This proof leads to the algorithm for construction of an
asymptotically optimal sequence of triangulations. Section 5
contains the proof of the lower estimate in Theorem~1.

\section{Preliminaries.} \label{3.1}
In order to investigate the asymptotic behavior of the optimal error
of piecewise linear interpolation of an arbitrary
function from the class $C^2(D)$ 
we shall use linear interpolation of the piecewise quadratic functions which
appear as an intermediate approximations of $f$.

Some of the facts we shall present in this section are quite easy to
see. However, we shall prove some of them, first of all for
completeness, and secondly, because we shall use them in the
construction of an asymptotically optimal sequence of
triangulations.

Let us define the modulus of continuity of $g \in C(D)$ by
$$
\omega(g, \delta):=\sup \{|g(x,y)-g(x',y')|: \;\; |x-x'|\le
\delta;\; |y-y'| \le \delta;\; (x,y), (x',y') \in D \}.
$$
For the function $f\in C^2(D)$ set
\begin{equation} \label{m1}
\omega(\delta):=\max\{\omega(f_{xx},\delta),\, \omega(f_{yy},\delta),\, \omega(f_{xy},\delta)  \}.
\end{equation}

\begin{lemma} \label{Lemma_Tay}
Let $f \in C^2(D)$. If $P_2=P_2(f;x,y;x_0,y_0)$ denotes the second
degree Taylor polynomial of $f$ at the point $(x_0,y_0)$ inside the
square $D_h \subset D$ with side length equal to $h$, then we have
the following estimate:
$$
  \|f-P_2\|_{L_{\infty}(D_h)}\le2h^2\,\omega(h),
$$
where $\omega(t)$ is defined in (\ref{m1}).
\end{lemma}

This simple lemma can be proved similarly to Lemma~1 from~\cite{us}.
The following statement is almost obvious.

\begin{lemma}\label{invar}
  For the given quadratic function
\begin{equation}\label{xxxyyy}
  Q(x,y)=Ax^2+By^2+2Cxy,
\end{equation}
an arbitrary triangle $T$, and any $c\in\RR^2$, the $L_p$-errors ($0<p\le\infty$)
of linear interpolation of $Q(x,y)$ on $T$, $c+T$, and a triangle $\widetilde{T}$
which is symmetric to $T$ with respect to the midpoint of any side of $T$, are equal, i.e.
$$
  d(Q,T,L_p)=d(Q,c+T,L_p)=d(Q,\widetilde{T},L_p).
$$
\end{lemma}

For an arbitrary linear transformation $S:\RR^2\to\RR^2$ denote by
${\rm det}\,S$ the determinant of the matrix of this transformation.

\begin{lemma} \label{invar1}
Consider a non-singular affine mapping $F=c+\widetilde{F}$, where
$c\in\R^2$ and $\widetilde{F}$ is a linear transformation. Then for
any quadratic function~(\ref{xxxyyy}) and any triangle $T$ we have
$|F^{-1}(T)|=\displaystyle\frac{|T|}{{\rm det}\,\widetilde{F}}$ and
$$
  d(Q\circ F,F^{-1}(T),L_p)=d(Q,T,L_p)\cdot\left(\frac{1}{{\rm det}\,\widetilde{F}}\right)^{\frac 1{p}}.
$$
\end{lemma}

This lemma can be easily verified by a routine change of variables.

Assume now that the function~(\ref{xxxyyy}) has a positive Hessian,
i.e. $AB-C^2>0$. Let us find the eigenvalues and unit eigenvectors
of the matrix of quadratic form ~(\ref{xxxyyy}).

For eigenvalues we have
$$   
\lambda_{\max}=\frac{A+B}{2}+\sqrt{\left(\frac{A+B}{2}\right)^2-(AB-C^2)}, \qquad
 \lambda_{\min}=\frac{A+B}{2}-\sqrt{\left(\frac{A+B}{2}\right)^2-(AB-C^2)}.
$$   
Observe that $0<\lambda_{\min}<\lambda_{\max}$.
In addition, note that
$$
  \lambda_{\min}\lambda_{\max}=AB-C^2.
$$
Let $(\xi_1, \xi_2)\in \SSs^1$ be an eigenvector of $Q(x,y)$
corresponding to the eigenvalue $\lambda_{\max}$. Then $(-\xi_2,
\xi_1) \in \SSs^1$ is an eigenvector corresponding to the eigenvalue
$\lambda_{\min}$.

\begin{lemma}\label{qua}
  For the quadratic form~(\ref{xxxyyy}) such that $AB-C^2>0$ it follows that
\begin{equation}\label{form}
  d(Q,T,L_p)\ge C_p^+|T|^{1+\frac 1p}\sqrt{\lambda_{\min}\lambda_{\max}}.
\end{equation}
\end{lemma}

{\bf Proof.} Recall that $\overline{Q}(x,y)=x^2+y^2$. Obviously,
$\overline{Q}$ is the canonical form of $Q$. We obtain this
canonical form by the following two linear transformations $F_1$ and
$F_2$
\begin{equation}\label{H1}
  F_1:\quad x=\xi_1 x'-\xi_2 y',\quad y=\xi_2 x'+\xi_1 y',
\end{equation}
and
\begin{equation}\label{H2}
  F_2:\quad x'=\frac{u}{\sqrt{\lambda_{\max}}},\quad y'=\frac v{\sqrt{\lambda_{\min}}}.
\end{equation}
Set $F:=F_1\circ F_2$. Note that
$$
  {\rm det}\,F=\frac{1}{\sqrt{\lambda_{\min}\lambda_{\max}}}
$$
and
$$
  Q\circ F(u,v)=u^2+v^2.
$$
Thus, in view of Lemma~\ref{invar1}
it follows that
$$
  d(Q,T,L_p)=d(Q\circ F,F^{-1}(T),L_p)\cdot({\rm det}\,F)^{\frac 1{p}}.
$$
Therefore, by the definition of the constant $C_p^+$ we obtain
$$
\begin{array}{rcl}
  \displaystyle\frac{d(Q,T,L_p)}{|T|^{1+\frac 1p}} & = & \displaystyle\frac{d(Q\circ F,F^{-1}(T),L_p)\cdot({\rm det}\,F)^{\frac 1p}}{|F^{-1}(T)|^{1+\frac 1p}\cdot({\rm det}\,F)^{1+\frac 1p}} \\[10pt]
  &  \ge & \displaystyle\inf\limits_{\widetilde{T}} \frac {d(\overline{Q},T,L_p)}{|\widetilde{T}|^{1+\frac 1p}}\cdot\sqrt{\lambda_{\min}\lambda_{\max}} = C_p^+\sqrt{\lambda_{\min}\lambda_{\max}}.
\end{array}
$$
Hence,
$$
  d_{Q,T,p}\ge C_p^+|T|^{1+\frac 1p}\sqrt{\lambda_{\min}\lambda_{\max}},
$$
which completes the proof. $\square$

Note that in Section~3 it will be shown that the infimum in
~(\ref{Cp}) is achieved on equilateral triangles. Thus, to obtain
the triangle on which the inequality~(\ref{form}) becomes equality,
we should take an arbitrary equilateral triangle $\widetilde T$, and
then the triangle $F(\widetilde{T})$ will be optimal.

In addition, we shall need the following two lemmas.
\begin{lemma} \label{L4}
Let us consider the collection of quadratic forms of type
(\ref{xxxyyy}) which satisfy the following conditions:
$$   
  0<A\le A^+, \;\; 0<B\le B^+,\;\;\hbox{and}\;\; H=AB-C^2 \ge C^+,
$$   
where $A^+, B^+, C^+$ are some positive numbers. Then for any such form
$$   
  \lambda_{\min} \geq \frac{1}{2}(A^++B^+)-\sqrt{\left (\frac{1}{2}(A^++B^+)\right )^2-C^+}>0.
$$   

\end{lemma}


To prove this lemma observe that the function
$g(u,v)=u-\sqrt{u^2-v}$ ($u>0$, $0<v\le1$) is decreasing in $u$ and
increasing in $v$.

\begin{lemma} \label{L11}
For the collection of quadratic forms satisfying the assumptions of
Lemma \ref{L4}, the ratio of the diameter of the optimal triangle
to the square root of the area of this triangle is bounded by the constant non-depending on $A^+$, $B^+$ and $C^+$.
\end{lemma}

This lemma follows from Lemma~5.

The next two results will be used in Section~5.

\begin{lemma}\label{second}
  Let $f\in C^2(D)$; $H(f;x,y)\ge C^+>0$ for all $(x,y)\in D$. If $\bar{n}$ is an arbitrary unit vector in the
  plane, then
\begin{equation}\label{const}
  \left|\frac{\partial^2 f}{\partial\bar{n}^2}\right|\ge D^+:=\frac{C^+}{2}\cdot\min{\left\{\frac1{\|f_{xx}\|_{\infty}};\frac1{\|f_{yy}\|_{\infty}}\right\}}.
\end{equation}
\end{lemma}

{\bf Proof. } Let $\bar{n}=(u,v)$ be an arbitrary unit vector in the plane. Then for an arbitrary point $(x,y)\in D$
$$
  \frac{\partial^2 f}{\partial \bar{n}^2}(x,y)=f_{xx}(x,y)u^2+2f_{xy}(x,y)uv+f_{yy}(x,y)v^2.
$$
Note that functions $f_{xx}$ and $f_{yy}$ have the same sign on $D$.
Without loss of generality we may assume that $f_{xx}(x,y)>0$ for
all $(x,y)\in D$. Since $u^2+v^2=1$ then either $u^2\ge\frac 12$, or
$v^2\ge\frac 12$. If $u^2>\frac 12$ then
$$
  \begin{array}{rcl}
  \displaystyle\frac{\partial^2 f}{\partial \bar{n}^2}(x,y) & = & \displaystyle\left(\frac{f_{xy}^2(x,y)}{f_{yy}(x,y)}u^2 + 2f_{xy}(x,y)uv+f_{yy}(x,y)v^2\right) + \left(f_{xx}(x,y) - \frac{f_{xy}^2(x,y)}{f_{yy}(x,y)}\right)u^2 \\[10pt]
  & \ge & \displaystyle\left(f_{xx}(x,y) - \frac{f_{xy}^2(x,y)}{f_{yy}(x,y)}\right)u^2 \ge \frac{C^+}{2\|f_{yy}\|_{\infty}}\ge D^+.
  \end{array}
$$
Similarly, if $v^2\ge\frac12$
$$
  \frac{\partial^2 f}{\partial \bar{n}^2}(x,y)\ge \frac{C^+}{2\|f_{xx}\|_{\infty}}\ge D^+.
$$
Thus, we have obtained the desired inequality. $\square$

\begin{lemma}\label{third}
  Let $f\in C^2(D)$; $H(f;x,y)\ge C^+>0$ for all $(x,y)\in D$. Then for any triangle $T$
$$
  d(f,T,L_p)\ge \frac{D^+}2\cdot d(\overline{Q},T,L_p),
$$
where constant $D^+$ was defined in~(\ref{const}).
\end{lemma}

{\bf Proof. } Let $A$, $B$, and $C$ be the vertices of the triangle $T$.
Set
$$
  n(x,y):=l(f,T)(x,y)-f(x,y), \quad (x,y)\in T,
$$
and
$$
  m(x,y):=\frac{D^+}{2}\cdot\left(l(\overline{Q},T)(x,y)-\overline{Q}(x,y)\right), \quad (x,y)\in T.
$$
Obviously, $n(x,y)\ge0$ and $m(x,y)\ge 0$ for all $(x,y)\in T$.
Moreover,
$$
  n(A)=n(B)=n(C)=m(A)=m(B)=m(C)=0.
$$
Let us consider the function
$$
  g(x,y):=n(x,y)-m(x,y), \quad (x,y)\in T.
$$
For the Hessian of $g$ we have
$$
\begin{array}{rcl}
  H(g;x,y) & = & \displaystyle\left(-f_{xx}(x,y)+D^+\right)\left(-f_{yy}(x,y)+D^+)-(f_{xy}^2(x,y)\right) \\[10pt]
  & = & \displaystyle f_{xx}(x,y)f_{yy}(x,y)-f_{xy}^2(x,y)+(D^+)^2-D^+\cdot\left (f_{xx}(x,y)+f_{yy}(x,y)\right) \\[10pt]
  &\ge& \displaystyle C^+ + (D^+)^2-D^+\cdot\left(f_{xx}(x,y)+f_{yy}(x,y)\right)\ge (D^+)^2.
\end{array}
$$
Therefore, $g$ is concave on the triangle $T$, since
$g_{xx}(x,y)\le0$ for all $(x,y)\in T$. It follows that $g(x,y)\ge0$
for all $(x,y)\in T$. Thus, $n(x,y)\ge m(x,y)$ for all $(x,y)\in T$,
from which the desired inequality easily follows. $\square$

\section{The proof of Theorem~2 for $1\le p<\infty$.}

In this section we shall provide the proof of Theorem 2 in the case when $1\le p<\infty$  and show the dependence
of the error of interpolation of the quadratic function on a
triangle on the geometry of the triangle.

\begin{lemma}\label{fff} Let $1\le p<\infty$. Then for
every non-equilateral triangle $T$ there exists a triangle $\overline{T}$ such that
$$
  \frac{d(\overline{Q},T,L_p)}{|T|^{1+\frac{1}{p}}}>\frac{d(\overline{Q},\overline{T},L_p)}{|\overline{T}|^{1+\frac 1p}}.
$$
\end{lemma}

{\bf Proof. } Assume that $T$ is not equilateral, i.e. $|AB|\ne
|BC|$. Let $S$ be an arbitrary rotation of the plane. Then,
obviously, $\overline{Q}\circ S\equiv \overline{Q}$. We may assume
that $A=(-1,\,0)$, $B=(a,\,b)$, and $C=(1,\,0)$, where $a\in\RR$ and
$b>0$.

Let $M=(x_M,y_M)$ be the center of the circle circumscribing
triangle $T$, and let $R$ be its radius. Clearly, the point $M$ has
the following coordinates:
$$
  M=\left(0,\,\frac{b^2+a^2-1}{2b}\right),
$$
and for the radius $R$ we have:
$$
  R^2=\left(\frac{b^2+a^2-1}{2b}\right)^2+1.
$$
Then for the error we obtain
\begin{equation}\label{1step}
  \begin{array}{rcl}
  d^p(\overline{Q},T,L_p) & = & d^p(\overline{Q}(\cdot,\cdot-y_M),T,L_p) \\[10pt]
                          & = & \displaystyle\int\limits_T\left|\left(\frac{b^2+a^2-1}{2b}\right)^2+1-x^2-\left(y-\frac{b^2+a^2-1}{2b}\right)^2\right|^p\,dx\,dy.
  \end{array}
\end{equation}

Let $T'$ be the triangle with the vertices $A$, $B'=(-a,\,b)$, and $C$, and let $\overline{T}$
be the triangle with the vertices $A$, $\overline{B}=(0,\,b)$, and $C$.
Obviously,
\begin{equation}\label{area1}
  |T|=|T'|=|\overline{T}|.
\end{equation}

Let us show that $\overline{T}$ is the desired triangle. To this
end, we consider the linear transformation $F$, which is determined
by the matrix
$$
  \left(\begin{array}{cc} 1 & \frac {a}{b} \\ 0 & 1 \end{array}\right).
$$
Note that $F$ transforms the triangle $\overline{T}$ into the
triangle $T$. Hence, we obtain
\begin{equation}\label{2step}
\begin{array}{rcl}
  d^p(\overline{Q},T,L_p) & = & d^p(\overline{Q}\circ F,\overline{T},L_p) \\[10pt]
                          & = & \displaystyle\int\limits_{\overline{T}}\left|1-u^2-2\frac{a}{b}\cdot uv-\frac{a^2}{b^2}\cdot v^2-v^2+\frac{b^2+a^2-1}{b}\cdot v\right|^p\,du\,dv.
\end{array}
\end{equation}
Note that
$$
  d(\overline{Q},T,L_p)=d(\overline{Q},T',L_p).
$$
Therefore, in view of~(\ref{area1}) and the triangle inequality, we obtain
$$
 \begin{array}{rcl}
  \frac{d(\overline{Q},T,L_p)}{|T|^{1+\frac1p}} & = & \displaystyle\frac{d(\overline{Q},T,L_p)+d(\overline{Q},T',L_p)}{2|T|^{1+\frac 1p}} \\[10pt]
                                                & \ge & \displaystyle\frac{1}{|\overline{T}|^{1+\frac 1p}}\left(\int\limits_{\overline{T}}\left|1-u^2-\frac{a^2}{b^2}\cdot v^2-v^2+\frac{b^2+a^2-1}{b}\cdot v\right|^p\,du\,dv\right)^{\frac 1p}.
  \end{array}
$$
Since $(u,v)\in \overline{T}$, we have $v\le b$. Using
~(\ref{1step}) with $a=0$ and equality~(\ref{2step}), it follows
that
$$
  \frac{d(\overline{Q},T,L_p)}{|T|^{1+\frac 1p}} > \frac{1}{|\overline{T}|^{1+\frac 1p}}\left(\int\limits_{\overline{T}}\left|1-u^2-v^2+\frac{b^2-1}{b}\cdot v\right|^p\,du\,dv\right)^{\frac 1p}
  =\frac{d(\overline{Q},\overline{T},L_p)}{|\overline{T}|^{1+\frac
  1p}}. \qquad \square
$$

\begin{lemma}\label{frf}
  Let $T$ be an arbitrary triangle. Then for any $0<p<\infty$
$$
  d(\overline{Q},T,L_p)\ge\frac{{\rm diam}\,T}{2^{5+\frac{1}{p}}h(T)}\cdot|T|^{1+\frac 1p},
$$
where $h(T)$ denotes the minimal height of the triangle $T$.
\end{lemma}

{\bf Proof. } Let $A$, $B$, $C$ be the vertices of the triangle $T$,
and let $|AB|={\rm diam}\,T$. We may assume that $A=(-a,0)$,
$B=(a,0)$ and $C=\left(b,h(T)\right)$, where $2a=|AB|$ and
$b\in(-a,a)$. Consequently, $a\cdot h(T)=1$. Without loss of
generality we may assume that $b\ge 0$. Let us consider the
trapezium $PQST$ with vertices $P=\left(\frac{b-3a}{4},\frac
1{4a}\right)$, $Q=\left(\frac{3b-a}4,\frac3{4a}\right)$,
$S=\left(\frac a2,0\right)$, and $T=\left(-\frac a2,0\right)$. Note
that the area of the trapezium $PQRS$ is equal to $\frac 12 |T|$.

Let $M=(x_M,y_M)$ be the center of the circumscribed circle of the
triangle $\overline{T}$, and let $R$ be its radius. Obviously,
$$
  M=\left(0,\,\frac{a}{2}\cdot(b^2+\frac{1}{a^2}-a^2)\right)\quad{\rm and}\quad R^2=a^2+y_M^2.
$$
Then
$$
  L:=\min{\{R^2-|MS|^2;\,R^2-|MT|^2;\,R^2-|MP|^2;\,R^2-|MQ|^2\}}\ge\frac{a^2}{16}.
$$
Recall that
$$
  d^p(\overline{Q},T,L_p)=\int\limits_{T}\left|R^2-(x-x_M)^2-\left(y-y_M\right)^2\right|^p\,dx\,dy.
$$
Moreover,
$$
  d^p(\overline{Q},T,L_p)\ge\int\limits_{PQST}\left|R^2-(x-x_M)^2-\left(y-y_M\right)^2\right|^p\,dx\,dy
  \ge L^p\int\limits_{PQST}dx\,dy\ge\frac{a^{2p}\cdot|T|}{2\cdot16^p}.
$$
Hence,
$$
  d^p(\overline{Q},T,L_p)\ge\frac{|AB|^{2p}|T|}{2^{6p+1}}=\frac{|AB|^{p}|T|^{p+1}}{2^{5p+1}h^p(T)}.\quad \square
$$

In view of Lemma~\ref{frf} we can derive that the infimum in the
right hand side of~(\ref{Cp}) exists. Moreover, by Lemma~\ref{fff},
this infimum is attained on the equilateral triangles only.

The proof of Theorem~2 in the case $p\in(0,1)$ can be found in the Appendix.

\section{Error of interpolation of $C^2$ functions by linear splines: upper estimate.} \label{3.2.1}


In this section we shall show that
$$
  \limsup_{N\to\infty}N\cdot R_N(f,L_{p,\Omega})\le\frac {C^+_p}{2}  \left ( \displaystyle \int_{D}H(f;x,y)^{\frac{p}{2(p+1)}} \Omega(x,y)^{\frac{p}{p+1}} dxdy\right)^{\frac{p+1}{p}}.
$$
In order to do so we are going to construct a sequence of
triangulations $\{\triangle_N^*\}_{N=1}^{\infty}$ such that
$$
  \limsup_{N\to\infty}N\cdot\|f-s(f,\triangle_N^*)\|_{L_{p,\Omega}}\le\frac {C^+_p}{2}  \left ( \displaystyle \int_{D}H(f;x,y)^{\frac{p}{2(p+1)}} \Omega(x,y)^{\frac{p}{p+1}} dxdy\right)^{\frac{p+1}{p}}.
$$

For a fixed $\varepsilon \in (0,1)$ and for every
$N \in \NN$ we define
\begin{equation} \label{m_N}
m_N=m_N(\varepsilon):=\min\left \{ m>0: \;\; \frac{2}{m^2}\, \omega\left (\frac{1}{m}\right ) \le \frac{\varepsilon}{N}
\right \},
\end{equation}
where $\omega(\delta)$ is the function defined in (\ref{m1}).

Observe that $m_N \to \infty$ as $N \to \infty$. In addition, note that
\begin{equation} \label{star}
\frac{N}{m_N^2} \to \infty, \;\;\; N \to \infty,
\end{equation}
i.e. $m_N=o\left(\sqrt{N}\right)$ as $N\to \infty$ (see, for
instance, Section~4 from~\cite{us}).

Let us divide the square $D$ into squares with the side length equal
to $\frac{1}{m_N}$ that have sides parallel to the sides of $D$, and
denote the resulting squares by $D_i^N:=D_i^N(\varepsilon)$,
$i=1,\dots,m_N^2$, enumerated in an arbitrary order.

Let $(x_i^N,y_i^N)$, $i=1,\dots,m_N^2$, be the center of the square $D_i^N$.
In addition, for every $i=1,\dots,m_N^2$ set
$$
  A_i^N:=\frac 12 f_{xx}(x_i^N, y_i^N),\;\; B_i^N:=\frac 12 f_{yy}(x_i^N,y_i^N),\;\;C_i^N:=f_{xy}(x_i^N, y_i^N)
$$
and
$$
  \overline{\Omega}_i^N:=\sup_{(x,y)\in D_i^N} \Omega(x,y).
$$
Note that
$$
  H(x_i^N,y_i^N):=H(f;x_i^N,y_i^N)=4(A^N_iB_i^N-(C_i^N)^2), \quad i=1,\dots,m_N^2.
$$
Set
$$
  Q_i^N(x,y):=A_i^N x^2+2C_i^N xy + B_i^N y^2,\quad (x,y)\in\RR^2.
$$

Now for the fixed $\varepsilon$ and for all $N$ large enough we will
construct an appropriate triangulation $\triangle_N(\varepsilon)$ of
$D$ consisting of $N$ triangles. To this end, we shall construct the
triangulations $\triangle_N^i:=\triangle_N^i(\varepsilon)$ of
squares $D_i^N$ depending on the eigenvalues of the quadratic forms
$Q_i^N$, $i=1,\dots,m_N^2$. After this we shall ``glue'' these
triangulations to obtain the triangulation
$\triangle_N(\varepsilon)$ of $D$.

Everywhere below $k_1, k_2, \dots,$ stand for constants
independent of $N$ and $\varepsilon$.

Note that for an arbitrary triangle $T\in D$ and for an arbitrary continuous function $g\in C(D)$
\begin{equation}\label{nerr}
  d(g,T,L_{p,\Omega})\le\|\Omega\|_{\infty}\|g\|_{L_{\infty}(T)}|T|^{\frac 1p}.
\end{equation}

For every $i=1,\ldots,m_N^2$ let
$$
  f_{N,i}(x,y)=P_2(f;x,y;x_i^N,y_i^N),
$$
where $P_2(f;x,y;a,b)$ is the Taylor polynomial of degree 2 of $f$
constructed at the point $(a,b)$.

For every $i=1,\dots,m_N^2$ set
$$
n_i^N=n_i^N(\varepsilon):=\left [ \frac{N(1-\varepsilon){H(x_i^N,
y_i^N)}^{\frac{p}{2(p+1)}}(\overline{\Omega}_i^N)^{\frac{p}{p+1}}}
{\displaystyle \sum_{j=1}^{m_N^2} {H(x_j^N,
y_j^N)}^{\frac{p}{2(p+1)}}(\overline{\Omega}_j^N)^{\frac{p}{p+1}}}\right
] + 1.
$$

Observe that all $n_i^N \to \infty$ when $N \to \infty$. This
follows from the obvious estimate
\begin{equation} \label{ni_est}
n_i^N \geq \left [ \frac{N(1-\varepsilon)
(C^+)^{\frac{p}{2(p+1)}} \displaystyle \min_{(x,y) \in D} \{
\Omega(x,y) \}^{\frac{p}{p+1}}}
{m_N^2\|H\|_{\infty}^{\frac{p}{2(p+1)}}\|\Omega\|_{\infty}^{\frac{p}{p+1}}}\right]+1,
\end{equation}
combined with~(\ref{star}), and $\min\limits_{(x,y) \in D}{
\Omega(x,y)}>0$.

For an arbitrary triangle $T\subset \RR^2$ let us consider the
triangle $\widetilde{T}$ which is symmetric to $T$ with respect to
the midpoint of one of the sides of $T$. Therefore,
$\Pi:=T\cup\widetilde{T}$ is a parallelogram. By tilling the plane
with the help of $\Pi$ we shall obtain a triangulation of the plane
$G(T)$.

Given $n_i^N$ for each square $D_i^N$, $i=1,\dots,m_N^2$, we
construct the triangulations $\triangle_N^i$ of the square $D_i^N$
as follows:

\begin{enumerate}
\item We consider transformations $F_{1}$ and $F_{2}$
of the form  (\ref{H1}) and (\ref{H2}) respectively, corresponding
to the quadratic function $Q_i^N(x,y)$. Set $F:=F_1\circ F_2$.

\item We take an arbitrary triangle $T$ which solves problem (\ref{Cp}), and consider the triangle $F(T)$.

\item Next we define $T_i^N$ to be a re-scaling of $F(T)$ so that
$$
|T_i^N|= \frac{1}{m_N^2n_i^N}.
$$

\item With the help of the triangle $T_i^N$ we generate the triangulation $G(T_i^N)$ of the plane as described above.

\item Triangles from $G(T_i^N)$ which lie completely in the square $D_i^N$ we include into triangulation $\triangle_N^i$.
Let us consider the triangles from $G(T_i^N)$ having common points
with the boundary of $D_i^N$. The intersection of an arbitrary such
triangle with $D_i^N$ is a polygon with at most $7$ vertices. We
triangulate this polygon into at most $5$ triangles without adding
new vertices. All triangles which we have obtained are included into
triangulation $\triangle_N^i$.


\end{enumerate}

 Now let us ``glue'' triangulations $\triangle_N^i$. We shall do it
according to the following algorithm.
\begin{enumerate}
\item For every $i=1,\dots,m_N^2$ let $W_i^N$ be the set of the vertices of
triangulation $\triangle_N^i$ which lie on the boundary of $D_i^N$.

\item For arbitrary $i,j=1,\ldots,m_N^2$, $i\ne j$, set $S_{i,j}=D_i^N\cap D_j^N$.

\item Let $S_{i,j}\ne\emptyset$. For every triangle $T\in\triangle_N^i$
if $T\cap S_{i,j}$ is a non-empty segment then we subdivide this
triangle by joining the vertices of $T$ and $W_j^N\cap T$. Then we
include all these triangles in the triangulation
$\triangle_N(\varepsilon)$.

\end{enumerate}

Therefore, for every $\varepsilon\in(0,1)$ there exists
$N(\varepsilon)\in\NN$ such that for all $N>N(\varepsilon)$ we
present a construction of the triangulation
$\triangle_N(\varepsilon)$.

Obviously, for every $i=1,\dots,m_N^2$ the conditions of Lemma~6 are
satisfied for the quadratic form $Q_i^N$. Thus there exists a
constant $c_1$, independent of $N$, such that for every
$i=1,\dots,m_N^2$ and for any triangle $T$ from the triangulation
$\triangle_N^i$ we have
\begin{equation}\label{trtr}
  {\rm diam}\,T \le c_1 \sqrt{\frac{1}{m_N^2n_i^N}}
\end{equation}
for all $N$ large enough.

For all $i=1,\dots,m_N^2$ denote by $\hat{n}_i^N$ the number of
triangles from $\triangle_N(\varepsilon)$ which lie in $D_i^N$. In
addition, we denote by $N_1=N_1(\varepsilon)$ the number of
triangles in the triangulation $\triangle_N(\varepsilon)$, and let
$\overline{N}_1=\overline{N}_1(\varepsilon)$ be the number of
triangles from the triangulation $\triangle_N(\varepsilon)$ which
have nonempty intersection with
$\bigcup\limits_{i=1}^{m_N^2}\partial D_i^N$.

From~(\ref{trtr}) it follows that for all $N$ large enough (without
loss of generality we may assume that this is true for all
$N>N(\varepsilon)$)
\begin{equation}\label{asd}
  \hat{n}_i^N\le(1+\varepsilon)n_i^N   \quad{\rm for\;all}\quad i={1,\dots,m_N^2},
\end{equation}
$$
  N_1\le (1-\varepsilon)N\quad{\rm and}\quad \overline{N}_1\le \varepsilon N.
$$

In view of construction of the triangulation
$\triangle_N(\varepsilon)$, for every triangle $T \in
\triangle_N(\varepsilon)$ which lies completely in ${\rm
int}\,D_i^N$ (the interior of the set $D_i^N$) we have
\begin{equation}\label{lll}
  d^p(f_{N,i},T,L_{p}) = d^p(f_{N,i},T_i^N,L_{p}) = \left(\frac{C^+_p}{2}\right)^p \frac{H(x_i^N,y_i^N)^{p/2}}{(m_N^2  n_i^N)^{p+1}}.
\end{equation}
Then for the error of interpolating $f_{N,i}$ on the square $D_i^N$ at the vertices of $\triangle_N^i$
we obtain that
$$
  \|f_{N,i}-s(f_{N,i},\triangle_N^i)\|_{L_{p,\Omega}(D_i^N)}^p\le (\overline{\Omega}_i^N)^p\sum_{{T\in \triangle_N(\varepsilon),\atop T\subset {\rm int}\,D_i^N} } d^p(f_{N,i},T,L_{p})\le \hat{n}_i^N (\overline{\Omega}_i^N)^p d^p(f_{N,i},T_i^N,L_{p}).
$$
Hence, in view of~(\ref{asd}) and~(\ref{lll}),
$$
V_1:=\sum_{i=1}^{m_N^2}\|f_{N,i}-s(f_{N,i}, \triangle^i_N)\|_{L_{p,\Omega}(D_i^N)}^p\le (1+\varepsilon)\left(\frac{C^+_p}{2}\right)^p \displaystyle \sum_{i=1}^{m_N^2}n_i^N \frac{(\overline{\Omega}_i^N)^p {H(x_i^N,y_i^N)}^{p/2}}{(m_N^2n_i^N)^{p+1}}.
$$
Using the definition of $n_i^N$ we obtain {\small
$$
\begin{array}{rcl}
V_1 & \le & \displaystyle\left(\frac{C^+_p}{2}\right)^p\frac{1+\varepsilon}{m_N^{2p+2}}
\displaystyle \sum_{i=1}^{m_N^2} {H(x_i^N, y_i^N)}^{p/2}
(\overline{\Omega}_i^N)^p\left ( \frac {\displaystyle
\sum_{j=1}^{m_N^2} {H(x_j^N,
y_j^N)}^{\frac{p}{2(p+1)}}(\overline{\Omega}_j^N)^{\frac{p}{p+1}}}{
N(1-\varepsilon){H(x_i^N,
y_i^N)}^{\frac{p}{2(p+1)}}(\overline{\Omega}_i^N)^{\frac{p}{p+1}}}\right
)^{p} \\ [10pt]
 & = & \displaystyle\left(\frac{C^+_p}{2}\right)^p\frac{1+\varepsilon}{(1-\varepsilon)^p N^p
m_N^{2(p+1)}} \left ( \displaystyle \sum_{j=1}^{m_N^2} H(x_j^N,
y_j^N)^{\frac{p}{2(p+1)}}(\overline{\Omega}_j^N)^{\frac{p}{p+1}}\right
)^{p+1}\\[10pt]
& = & \displaystyle\left(\frac{C^+_p}{2}\right)^p\frac{1+\varepsilon}{(1-\varepsilon)^p N^p} \left ( \displaystyle \frac{1}{m_N^2}\sum_{j=1}^{m_N^2} H(x_j^N,
y_j^N)^{\frac{p}{2(p+1)}}(\overline{\Omega}_j^N)^{\frac{p}{p+1}}\right
)^{p+1}.
\end{array}
$$}
Hence, by Riemann integrability of functions $H(f;x,y)$ and
$\Omega(x,y)$
$$
\begin{array}{rcl}
\displaystyle\lim_{N\to \infty}\frac{1}{m_N^2}\displaystyle
\sum_{j=1}^{m_N^2}H(x_j^N,
y_j^N)^{\frac{p}{2(p+1)}}(\overline{\Omega}_j^N)^{\frac{p}{p+1}} & = & \displaystyle \lim_{N\to \infty}\sum_{j=1}^{m_N^2}|D_j^N|H(x_j^N, y_j^N)^{\frac{p}{2(p+1)}}(\overline{\Omega}_j^N)^{\frac{p}{p+1}} \\[10 pt]
& = &
\displaystyle\int_{D}H(f;x,y)^{\frac{p}{2(p+1)}}\Omega(x,y)^{\frac{p}{p+1}}dxdy.
\end{array}
$$
Hence, for all $N$ large enough we obtain
$$
V_1 \le \left(\frac{C^+_p}{2N}\right)^p(1+k_1\varepsilon)\left (
\int_{D}H(f;x,y)^{\frac{p}{2(p+1)}}\Omega(x,y)^{\frac{p}{p+1}}dxdy\right )^{p+1}.
$$
Note that for all $\varepsilon\in(0,1)$ and for all $N>N(\varepsilon)$
we have
\begin{equation}\label{l2}
  R_N(f,L_{p,\Omega})\le \|f-s(f,\triangle_N(\varepsilon))\|_{p,\Omega}.
\end{equation}

Let us consider the case $1\le p<\infty$.
Observe that
\begin{equation}\label{wlw}
\begin{array}{rcl}
  \|f-s(f,\triangle_N(\varepsilon))\|_{p,\Omega}^p  & = &  \displaystyle\sum_{i=1}^{m_N^2}\|f-s(f,\triangle_N\varepsilon))\|^p_{p,\Omega} \le \displaystyle\sum_{i=1}^{m_N^2}\Big[\|f-f_{N,i}\|_{p,\Omega}\\[10pt]
  & +& \displaystyle\|f_{N,i}-s(f_{N,i}, \triangle_N(\varepsilon))\|_{p,\Omega} +\|s(f_{N,i},\triangle_N(\varepsilon))-s(f, \triangle_N(\varepsilon))\|_{p,\Omega}\Big]^p.
\end{array}
\end{equation}
Obviously, for every $i=1,\ldots,m_N^2$
$$
  \|s(f_{N,i},\triangle_N(\varepsilon))-s(f, \triangle_N(\varepsilon))\|_{p,\Omega}\le \frac 1{m_N^2}\|f-f_{N,i}\|_{\infty,\Omega}.
$$
and
$$
  \|f-f_{N,i}\|_{p,\Omega}\le \frac 1{m_N^2} \|f-f_{N,i}\|_{\infty,\Omega}.
$$
In addition, by Lemma~\ref{Lemma_Tay} and the definition of $m_N$ we
have
$$
  \|f-f_{N,i}\|_{\infty,\Omega} \le \frac{2\|\Omega\|_{\infty}}{ m_N^2}\, \omega\left(\frac{1}{m_N}\right) \le \frac{\varepsilon}{N}\|\Omega\|_{\infty}.
$$
Thus,
$$
  \|f-s(f,\triangle_N(\varepsilon))\|_{p,\Omega}^p\le\displaystyle\sum_{i=1}^{m_N^2}\Big[\|f_{N,i}-s(f_{N,i}, \triangle_N(\varepsilon))\|_{p,\Omega} + \frac{2\varepsilon}{Nm_N^2}\|\Omega\|_{\infty}\Big]^p.
$$
Hence,
$$
\begin{array}{rcl}
  \|f-s(f,\triangle_N(\varepsilon))\|_{p,\Omega}^p & \le & \displaystyle\sum_{i=1}^{m_N^2}\|f_{N,i}-s(f_{N,i}, \triangle_N(\varepsilon))\|_{p,\Omega}^p + \frac{k_2\varepsilon}{N^p}. \\[10 pt]
        & \le & \displaystyle (1+k_3\varepsilon)\sum_{i=1}^{m_N^2}\|f_{N,i}-s(f_{N,i}, \triangle_N(\varepsilon))\|^p_{p,\Omega} =: (1+k_3\varepsilon) V_1.
\end{array}
$$
Therefore, there exists $\widetilde{N}(\varepsilon)\ge N(\varepsilon)$ such that for all
$N\ge \widetilde{N}(\varepsilon)$
\begin{equation}\label{vvv}
  \|f-s(f,\triangle_N(\varepsilon))\|_{p,\Omega}^p \le \left(\frac{C^+_p}{2N}\right)^p(1+k_4\varepsilon)\left (
\int_{D}H(f;x,y)^{\frac{p}{2(p+1)}}\Omega(x,y)^{\frac{p}{p+1}}dxdy\right )^{p+1}.
\end{equation}

Note that in the case $p\in(0,1)$ instead of ~(\ref{wlw}) we should
use the following implication of ~(\ref{norm_inverse}):
$$
  \begin{array}{rcl}
    \|f-s(f,\triangle_N(\varepsilon))\|_{p,\Omega}^p &  =  & \displaystyle\sum_{i=1}^{m_N^2}\|f-s(f,\triangle_N(\varepsilon))\|^p_{p,\Omega} \le \displaystyle\sum_{i=1}^{m_N^2}\|f-f_{N,i}\|^p_{p,\Omega} \\ [10 pt]
          &     &  + \|f_{N,i}-s(f_{N,i}, \triangle_N(\varepsilon))\|^p_{p,\Omega} +\|s(f_{N,i},\triangle_N(\varepsilon))-s(f, \triangle_N(\varepsilon))\|^p_{p,\Omega}.
  \end{array}
$$
The rest of the proof is analogous to the case $1\le p < \infty$.

Let $1>\varepsilon_{1}>\varepsilon_2>\dots$ be a decreasing sequence
of positive numbers which tends to zero as $k\to \infty$. Without
loss of generality we may assume that
$\{\widetilde{N}(\varepsilon_k)\}_{k=1}^{\infty}$ is an increasing
sequence. Then set
$$
  \triangle_N^*:=\triangle_N(\varepsilon_k)\quad {\rm if }\quad \widetilde{N}(\varepsilon_k)<N\le \widetilde{N}(\varepsilon_{k+1}), \quad k\in \NN,
$$
and let $\triangle_N^*$ be an arbitrary triangulation of $D$ if $1\le N\le \widetilde{N}(\varepsilon_1)$.

Therefore, in view of the definition of $R_N(f,L_{p,\Omega})$ and
inequalities~(\ref{l2}) and~(\ref{vvv}), for all $0<p<\infty$ and
for every $\widetilde{N}(\varepsilon_k)<N\le
\widetilde{N}(\varepsilon_{k+1})$ we have

$$
  R_N(f,L_{p,\Omega})\le \|f-s(f,\triangle_N^*)\|_{p,\Omega}
  \le  \frac{C^+_p}{2N}(1+k_5\varepsilon_k)\left( \int_{D}H(f;x,y)^{\frac{p}{2(p+1)}}\Omega(x,y)^{\frac{p}{p+1}}dxdy\right )^{\frac{p+1}{p}}.
$$
Since $\varepsilon_k\to0$, as $k\to\infty$, we obtain the desired
estimate. $\square$

\section{Error of interpolation of $C^2$ functions by linear splines: lower estimate.}

For an arbitrary triangle $T$ in the plane denote by $h(T)$, ${\rm
diam}\,T$, $U(T)$ and $|T|$ the minimal height of $T$, the length of
the longest side of $T$, an arbitrary vertex of the longest side of
$T$, and the area of $T$, respectively.

To prove the lower estimate we need to show that
\begin{equation}\label{formula1}
  \liminf_{N\to\infty} N R_N(f,L_{p,\Omega}) \ge \frac {C^+_p}{2}  \left ( \displaystyle \int_{D}H(f;x,y)^{\frac{p}{2(p+1)}} \Omega(x,y)^{\frac{p}{p+1}} dxdy\right)^{\frac{p+1}{p}}.
\end{equation}

Let $\underline{\Omega}=\min\limits_{(x,y)\in D}\Omega(x,y)$, and let $\{\triangle_N\}_{N=1}^{\infty}$
be an arbitrary sequence of triangulations of $D$.

For any $\varepsilon>0$ let us consider the following sets:
$$
  I_N(\varepsilon):=\Big\{i\in\{1,\ldots,N\}\;:\;\frac{h(T_i^N)}{{\rm diam}\,T_i^N} < \frac{c}{\varepsilon}\omega({\rm diam}\,T^N_i)\Big\},
$$
where
$$
  c:=\displaystyle\frac{8\|\Omega\|_{\infty}}{\frac{C_p^+}{2}\underline{\Omega}\sqrt{H(f;U(T_i^N))}},
$$
and
$$
  J_N(\varepsilon):=\Big\{i\in\{1,\ldots,N\}\;:\;{\rm diam}\,T_i^N>\frac {\varepsilon}{\sqrt[4]{N}}\Big\}.
$$
Note that the sets $I_N(\varepsilon)$ and $J_N(\varepsilon)$ might
have a nonempty intersection.

First, assume that there exists $\varepsilon_0>0$ and a subsequence
$\{N_k\}_{k=1}^{\infty}$ of positive integers such that
\begin{equation}\label{1assump}
  \sum_{i\in I_{N_k}(\varepsilon_0)\setminus J_{N_k}(\varepsilon_0)}|T^{N_k}_i|\ge\varepsilon_0.
\end{equation}

Let us show that in this case
$\liminf\limits_{k\to\infty}N_k\|f-s(f,\triangle_{N_k})\|_{p,\Omega}=\infty$.

Indeed, in view of Lemma~\ref{second} there exists a constant $D^+$
such that the derivative of the function $f$ in an arbitrary unit
direction is bounded away from zero by $D^+$. Thus, by
Lemma~\ref{third} we obtain that
$$
  \begin{array}{rcl}
  \|f-s(f,\triangle_{N_k})\|^p_{p,\Omega} & \ge & \displaystyle\sum_{i\in I_{N_k}(\varepsilon_0)}d^p(f,T_i^{N_k},L_{p,\Omega}) \\[10pt]
  & \ge & \displaystyle(\underline{\Omega})^p \sum_{i\in I_{N_k}(\varepsilon_0)\setminus J_{N_k}(\varepsilon_0)}d^p(f,T_i^{N_k},L_p) \\[10pt]
  & \ge & \displaystyle\left(\frac{\underline{\Omega}D^+}{2}\right)^p \sum_{i\in I_{N_k}(\varepsilon_0)\setminus J_{N_k}(\varepsilon_0)}d^p(\overline{Q},T_i^{N_k},L_p).
  \end{array}
$$
Hence, by Lemma~10 and the definition of $I_{N_k}(\varepsilon_0)$ we
have
$$
\begin{array}{rcl}
  \|f-s(f,\triangle_{N_k})\|^p_{p,\Omega} & \ge &\displaystyle \left(\frac{\underline{\Omega}D^+}{2}\right)^p \frac{1}{2^{5p+1}}\sum_{i\in I_{N_k}(\varepsilon_0)\setminus J_{N_k}(\varepsilon_0)}|T_i^{N_k}|^{p+1}\cdot\left(\frac{{\rm diam}\,T_i^{N_k}}{h(T_i^{N_k})}\right)^p\\[10pt]
  &\ge &\displaystyle\left(\frac{\varepsilon_0\underline{\Omega}D^+}{2c}\right)^p\frac{1}{2^{5p+1}}\sum_{i\in I_{N_k}(\varepsilon_0)\setminus J_{N_k}(\varepsilon_0)}\frac{|T_i^{N_k}|^{p+1}}{\omega^p({\rm diam}\,T_i^{N_k})}.
\end{array}
$$
Let ${\rm card}\,(F)$ denote the number of elements in the finite
set $F$. Since ${\rm
diam}\,(T_i^N)\le\frac{\varepsilon_0}{\sqrt[4]{N}}$, applying the
Jensen inequality and ~(\ref{1assump}), we obtain that
$$
  \begin{array}{rcl}
    \|f-s(f,\triangle_{N_k})\|^p_{p,\Omega} & \ge & \displaystyle\left(\frac{\varepsilon_0\underline{\Omega}D^+}{2c\cdot\omega\left(\frac{\varepsilon_0}{\sqrt[4]{N_k}}\right)}\right)^p\frac{1}{2^{5p+1}} \sum\limits_{i\in I_{N_k}(\varepsilon_0)\setminus J_{N_k}(\varepsilon_0)}|T_i^{N_k}|^{p+1} \\[10pt]
   & \ge & \displaystyle\left(\frac{\varepsilon_0\underline{\Omega}D^+}{2c\cdot\omega\left(\frac{\varepsilon_0}{\sqrt[4]{N_k}}\right)}\right)^p \frac{\left(\sum\limits_{i\in I_{N_k}(\varepsilon_0)\setminus J_{N_k}(\varepsilon_0)}|T_i^{N_k}|\right)^{p+1}}{2^{5p+1}\Big[{\rm card}\,(I_{N_k}(\varepsilon_0)\setminus J_{N_k}(\varepsilon_0))\Big]^p} \\ [10pt]
   & \ge & \displaystyle\left(\frac{\varepsilon_0\underline{\Omega}D^+}{2c\cdot\omega\left(\frac{\varepsilon_0}{\sqrt[4]{N_k}}\right)}\right)^p \frac{\varepsilon_0^{p+1}}{2^{5p+1}{N_k}^p},
  \end{array}
$$
The last inequality implies that $\liminf\limits_{k\to\infty}N_k
\|f-s(f,\triangle_{N_k})\|^p_{p,\Omega}=\infty$.

Assume now that there exists a number $\varepsilon_0>0$ and a
subsequence $\{N_k\}_{k=1}^{\infty}$ of positive integers such that
$$
  \sum_{i\in J_{N_k}(\varepsilon_0)}|T^{N_k}_i|\ge\varepsilon_0.
$$
As in the previous case we are going to show that
$\liminf\limits_{k\to\infty}N_k\|f-s(f,\triangle_{N_k})\|_{p,\Omega}=\infty$.

Similarly to above we can derive that
$$
  \|f-s(f,\triangle_{N_k})\|^p_{p,\Omega}\ge \left(\frac{\underline{\Omega}D^+}{2}\right)^p \frac{1}{2^{5p+1}}\sum_{i\in J_{N_k}(\varepsilon_0)}|T_i^{N_k}|^{p+1}\left(\frac{{\rm diam}\,T_i^{N_k}}{h(T_i^{N_k})}\right)^p.
$$
Therefore, in view of $|T|=\frac 12 {\rm diam T} h(T)$ for an
arbitrary triangle $T$ and the definition of
$J_{N_k}(\varepsilon_0)$, we have
$$
  \begin{array}{rcl}
  \|f-s(f,\triangle_{N_k})\|^p_{p,\Omega} & \ge& \displaystyle\left(\frac{\underline{\Omega}D^+}{4}\right)^p \frac{1}{2^{5p+1}}\sum_{i\in J_{N_k}(\varepsilon_0)}|T_i^{N_k}|({\rm diam}\,T_i^{N_k})^{2p}\\[10pt]
  &\ge&\displaystyle\left(\frac{\underline{\Omega}D^+}{4N_k^{\frac 12}}\right)^p \frac{\varepsilon_0^{2p}}{2^{5p+1}}\sum_{i\in J_{N_k}(\varepsilon_0)}|T_i^{N_k}|\ge \frac{\varepsilon_0^{2p+1}}{2^{5p+1}}\left(\frac{\underline{\Omega}D^+}{4N_k^{\frac 12}}\right)^p.
  \end{array}
$$
This implies that
$\liminf\limits_{k\to\infty}N_k\cdot\|f-s(f,\triangle_{N_k})\|^p_{p,\Omega}=\infty$.

Hence, in order to find a sequence of triangulations
$\{\triangle_N\}_{N=1}^{\infty}$ which provides the lowest value of
$\liminf\limits_{N\to\infty}N\cdot\|f-s(f,\triangle_N)\|_{p,\Omega}$,
we should investigate only the sequences of triangulations with the
following property: for an arbitrary $\varepsilon>0$ there exists
$N(\varepsilon)$ such that for all $N>N(\varepsilon)$
$$
  \sum_{i\in I_{N}(\varepsilon)\cup J_{N}(\varepsilon)}|T^{N}_i|<\varepsilon.
$$
Let $\{\triangle_N\}_{N=1}^{\infty}$ be a sequence of triangulations
satisfying this condition. For every $i=1,\ldots,N$ let $f_{N,i}$ be
the second degree Taylor polynomial of $f$ at the point
$U(T_i^{N})$. Then in view of Lemma~1 we obtain that
$$
  \|f-f_{N,i}\|_{L_{\infty}(T_i^{N})}\le 2({\rm diam}\,T_i^N)^2\omega({\rm diam}\,T_i^N)
$$
for all $i=1,\ldots,N$.
Therefore,
$$
  \|f-f_{N,i}\|_{L_{p,\Omega}(T^{N}_i)}\le \|\Omega\|_{\infty}\|f-f_{N,i}\|_{L_{\infty}(T_i^{N})}|T_i^{N}|^{\frac 1p}\le 2\|\Omega\|_{\infty}({\rm diam}\,T_i^N)^2\omega({\rm diam}\,T_i^N)|T_i^{N}|^{\frac 1p},
$$
and
$$
  \|l(f,T_i^{N})-l(f_{N,i},T_i^{N})\|_{L_{p,\Omega}(T^{N}_i)}\le 2\|\Omega\|_{\infty}({\rm diam}\,T_i^N)^2\omega({\rm diam}\,T_i^N)|T_i^{N}|^{\frac 1p}.
$$
Consequently, for all $i\not\in I_{N}(\varepsilon)$
$$
  \begin{array}{rcl}
  \|f-f_{N,i}\|_{L_{p,\Omega}(T^{N}_i)}& + &\|l(f,T_i^{N})-l(f_{N,i},T_i^{N})\|_{L_{p,\Omega}(T^{N}_i)}\\[10pt]
  &\le & \displaystyle4\|\Omega\|_{\infty}({\rm diam}\,T_i^N)^2\omega({\rm diam}\,T_i^N)|T_i^{N}|^{\frac 1p}\\[10pt]
  &\le &\displaystyle4\|\Omega\|_{\infty}({\rm diam}\,T_i^N)^2|T_i^{N}|^{\frac 1p}\frac{\varepsilon\cdot h(T_i^{N})}{c\cdot {\rm diam}\,T_i^N} = \frac{8\varepsilon\cdot\|\Omega\|_{\infty}}{c}|T_i^{N}|^{1 + \frac 1p}\\[10pt]
  &\le&\displaystyle \varepsilon\cdot \frac{C_p^+}{2}\cdot\underline{\Omega}\sqrt{H(f;U(T_i^{N}))}\cdot|T_i^{N}|^{1+\frac 1p}.
  \end{array}
$$
Set $\underline{\Omega}_i^{N}:=\inf\limits_{(x,y)\in
T_i^{N}}\Omega(x,y)$, $i=1,\ldots,N$. Obviously,
\begin{equation}\label{fof}
  \|f-s(f,\triangle_{N})\|_{p,\Omega}^p = \displaystyle\sum_{i=1}^{N}d^p(f,T_i^N,L_{p,\Omega})
  \ge \sum_{i\not\in I_{N}(\varepsilon)\cup J_{N}(\varepsilon)}d^p(f,T_i^N,L_{p,\Omega}).
\end{equation}

Applying the triangle inequality and the defintion of
sets $I_N(\varepsilon)$ and $J_N(\varepsilon)$, for all $1\le
p<\infty$ and $i=1,\ldots,N$ we have
$$
  \begin{array}{rcl}
  d^p(f,T_i^N,L_{p,\Omega}) & \ge & \displaystyle\Big(d(f_{N,i},T_i^N,L_{p,\Omega})\\[10pt]
  & & \begin{array}{ll}
        - & \;\displaystyle\|f-f_{N,i}\|_{L_{p,\Omega}(T^{N}_i)}-\|l(f,T^{N}_i)-l(f_{N,i},T^{N}_i)\|_{L_{p,\Omega}(T^{N}_i)}\Big)^p
      \end{array}\\[10pt]
  & \ge & \Big(\frac{C_p^+}{2}\sqrt{H(f;U(T_i^{N}))}\cdot\underline{\Omega}_i^{N}|T_i^{N}|^{1+\frac 1p} \\ [10pt]
        & & \begin{array}{ll}
           -&\;\displaystyle 4\|\Omega\|_{\infty}({\rm diam}\,T_i^N)^2\omega({\rm diam}\,T_i^N)|T_i^{N}|^{\frac 1p}\Big)^p
            \end{array}\\[10pt]
  &\ge & \displaystyle\left((1-\varepsilon)\cdot\frac{C_p^+}{2}\right)^p\Big(\sqrt{H(f;U(T_i^{N}))}\cdot\underline{\Omega}_i^{N}\Big)^p|T_i^{N}|^{p+1}.
  \end{array}
$$
Similarly, for $p\in(0,1)$ and $i=1,\ldots,N$ we have
$$
  d^p(f,T_i^N,L_{p,\Omega}) \ge \displaystyle\left((1-k_6\varepsilon)\cdot\frac{C_p^+}{2}\right)^p\Big(\sqrt{H(f;U(T_i^{N}))}\cdot\underline{\Omega}_i^{N}\Big)^p|T_i^{N}|^{p+1}.
$$
Therefore, in view of~(\ref{fof}), for all $p\in(0,\infty)$
$$
  \begin{array}{rcl}
  \|f-s(f,\triangle_{N})\|_{p,\Omega}^p &\ge& \displaystyle\left((1-k_7\varepsilon)\cdot\frac{C_p^+}{2}\right)^p\sum_{i\not\in I_{N}(\varepsilon)\cup J_{N}(\varepsilon)}\Big(\sqrt{H(f;U(T_i^{N}))}\cdot\underline{\Omega}_i^{N}\Big)^p|T_i^{N}|^{p+1} \\[10pt]
& \ge & \displaystyle\left(\frac{(1-k_7\varepsilon)C_p^+}{2(N-|I_{N}(\varepsilon)\cup J_{N}(\varepsilon)|)}\right)^p\\[10pt]
  & & \begin{array}{ll}
      \times & \;\displaystyle\left(\sum_{i\not\in I_{N}(\varepsilon)\cup J_{N}(\varepsilon)}\Big[H(f;U(T_i^{N}))\Big]^{\frac{p}{2(p+1)}}(\underline{\Omega}_i^{N})^{\frac{p}{p+1}}|T_i^{N}|\right)^{p+1}
      \end{array}\\[10pt]
  & \ge & \displaystyle\left(\frac{(1-k_7\varepsilon)C_p^+}{2N}\right)^p\left(\sum_{i\not\in I_{N}(\varepsilon)\cup J_{N}(\varepsilon)}\Big[H(f;U(T_i^{N}))\Big]^{\frac{p}{2(p+1)}}(\underline{\Omega}_i^{N})^{\frac{p}{p+1}}|T_i^{N}|\right)^{p+1}.
  \end{array}
$$
Let us divide each triangle $T_i^{N}$, $i\in I_{N}(\varepsilon)\cup
J_{N}(\varepsilon)$, into $n_i^{N}$ triangles $T_{i,j}^{N}$,
$j=1,\ldots,n_i^{N}$, enumerated in an arbitrary order, such that
${\rm diam}\,(T_{i,j}^{N})\to 0$ as $N\to\infty$ for all
$j=1,\ldots,n_i^{N}$. For every $i\not\in I_{N}(\varepsilon)\cup
J_{N}(\varepsilon)$ set $n_i^{N}=1$ and $T_{i,1}^{N}=T_i^{N}$. In
addition, set $\underline{\Omega}_{i,j}^{N}:=\inf\limits_{(x,y)\in
T_{i,j}^{N}}\Omega(x,y)$ for all $j=1,\ldots,n_i^{N}$ and
$i=1,\ldots,N$. Note, that
$\bigcup\limits_{i=1}^{N}\bigcup\limits_{j=1}^{n_i^{N}}T_{i,j}^{N}=D$,
and for every $i$ and $j$ we have that ${\rm
diam}\,(T_{i,j}^{N})\to0$ as $N\to\infty$. Then
$$
  \begin{array}{rcl}
  \|f-s(f,\triangle_{N})\|_{p,\Omega}^p & \ge & \displaystyle\left(\frac{(1-k_7\varepsilon)C_p^+}{2N}\right)^p\left(\sum_{i=1}^{N}\sum_{j=1}^{n_i^{N}}\Big[H(f;U(T_{i,j}^{N}))\Big]^{\frac{p}{2(p+1)}}(\underline{\Omega}_{i,j}^{N})^{\frac{p}{p+1}}|T_{i,j}^{N}|\right. \\[10pt]
  & & \begin{array}{ll}
      - & \;\displaystyle\left.(\|\Omega\|_{\infty}\sqrt{\|H\|_{\infty}})^{\frac{p}{p+1}}\sum_{i\in I_{N}(\varepsilon)\cup J_{N}(\varepsilon)}|T^{N}_i|\right)^{p+1}
      \end{array} \\[10pt]
  & \ge & \displaystyle\left(\frac{C_p^+}{2(1+k_8\varepsilon)N}\right)^p\left(\int\limits_DH^{\frac{p}{2(p+1)}}(f;x,y)\Omega^{\frac{p}{p+1}}(x,y)\,dx\,dy\right)^{p+1}.
  \end{array}
$$
 Therefore,
$$
  \liminf_{N\to\infty}N \|f-s(f,\triangle_{N})\|_{p,\Omega}\ge \frac{C_p^+}{2(1+k_8\varepsilon)}\left(\int\limits_DH^{\frac{p}{2(p+1)}}(f;x,y)
  \Omega^{\frac{p}{p+1}}(x,y)\,dx\,dy\right)^{1+\frac{1}{p}}
$$
and since $\varepsilon$ is arbitrary we obtain the desired
inequality. $\square$

\section{Appendix}

In this section we shall show that for $0<p<1$
$$
  C_p^+:=\inf_{T}\frac{d(\overline{Q},T,L_p)}{|T|^{1+\frac 1p}} = \left(\frac{4}{3\sqrt{3}}\right)^{1+\frac 1p}\left[\frac{\pi}{p+1}-6\int\limits_{1/2}^1 x(1-x^2)^p\arccos{\frac 1{2x}}\,dx\right]^{\frac 1p}.
$$
Note that
\begin{equation}\label{inf}
  C_p^+=\inf_{T,\, |T|=1}{d(\overline{Q},T,L_p)}.
\end{equation}
By $M$ and $R$ we shall denote the center and radius of the circle
circumscribing triangle $T$, respectively. We may assume that the
point $M$ is at the origin. Then
$$
  d^p(\overline{Q},T,L_p)=\int\!\!\!\int\limits_{T}(R^2-x^2-y^2)^p\,dx\,dy.
$$
\begin{figure}[here]
   \centering
   \includegraphics[width=2 in]{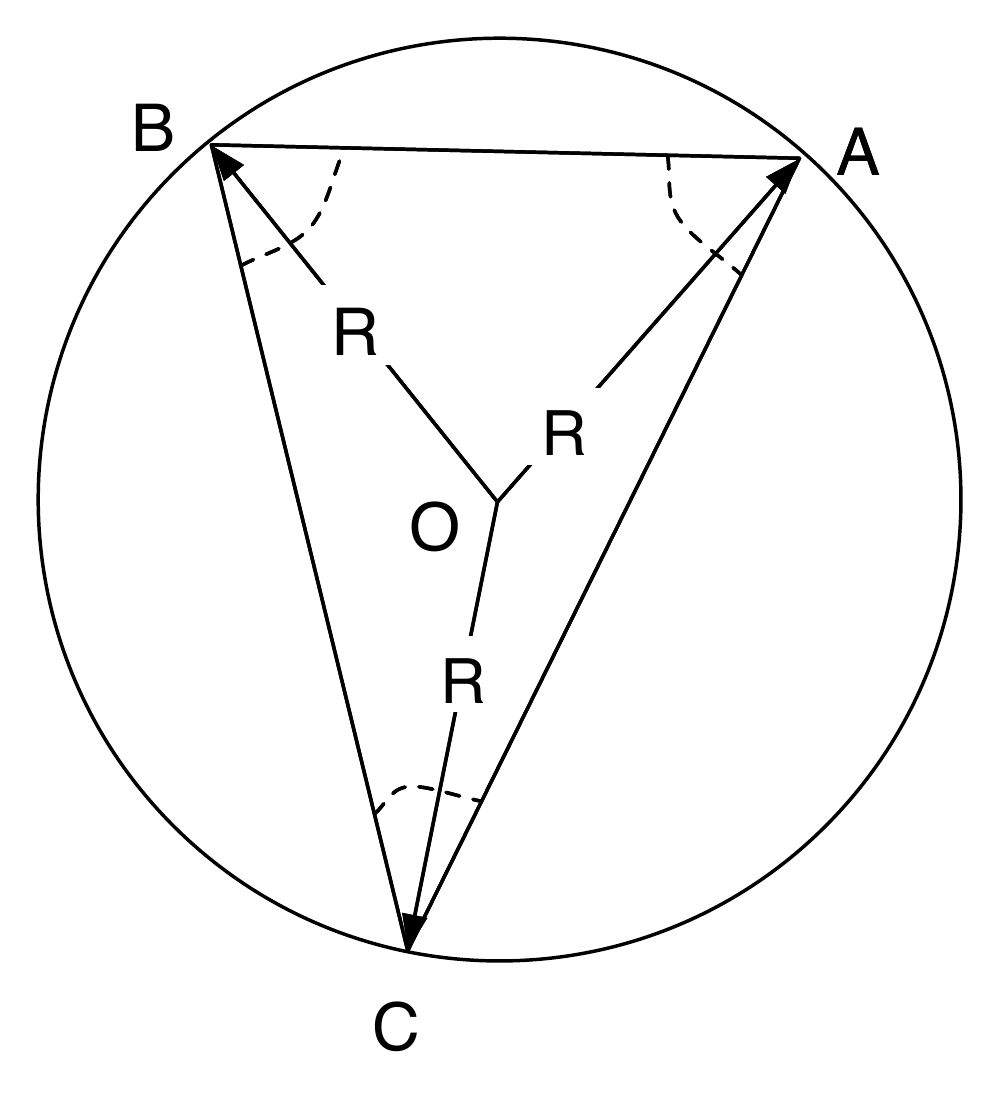}\qquad\qquad\qquad
   \includegraphics[width=2.2 in]{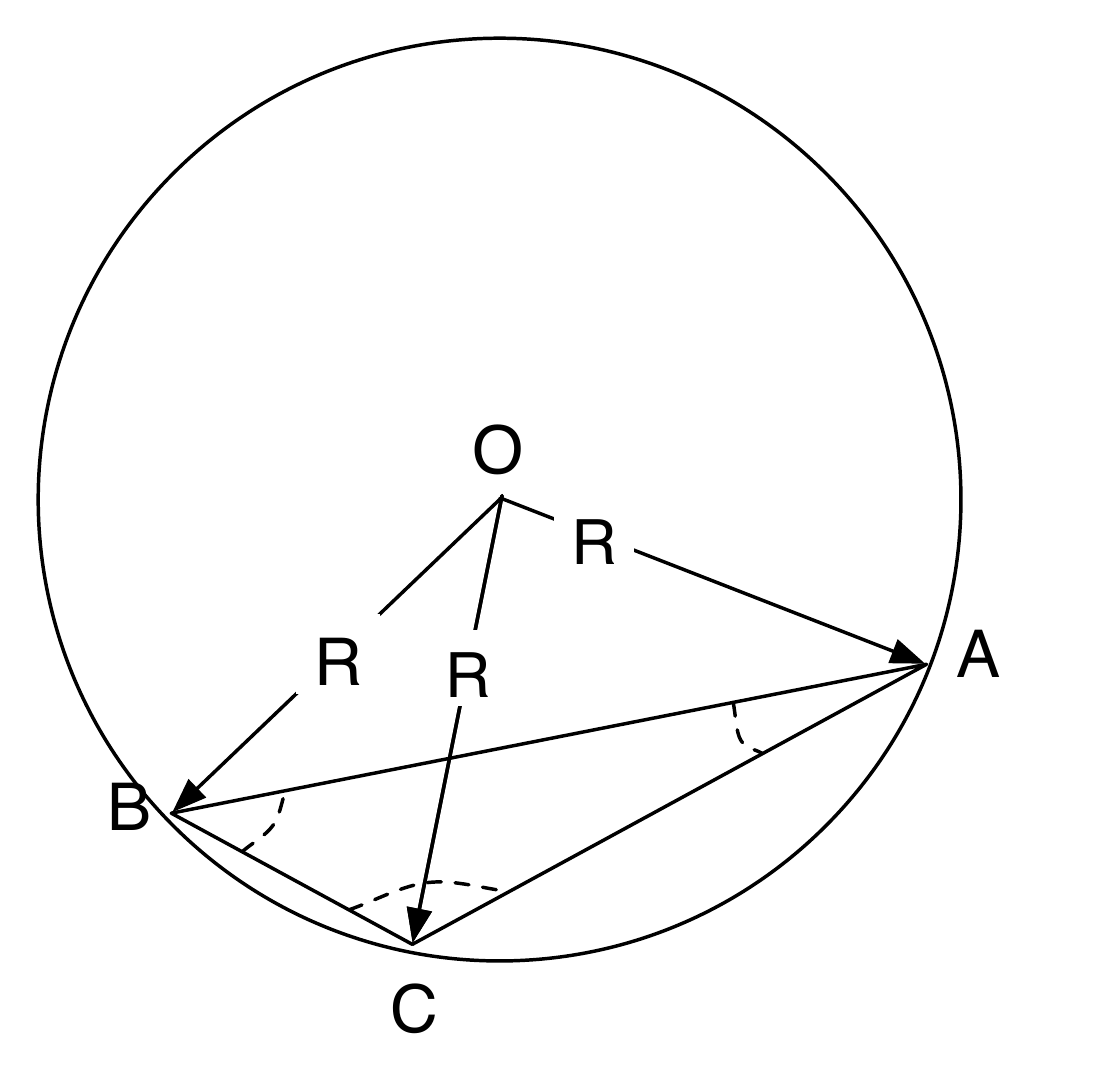}
   \caption{\it The triangles $T_A=\triangle BMC$, $T_B=\triangle CMA$, and $T_C=\triangle AMB$,
    in the cases of acute (left) and obtuse (right) triangle}
   \label{fig:triangle_1}
\end{figure}

Let $A$, $B$, $C$ be the vertices of the triangle $T$. By the same
letters we shall denote the angles corresponding to the vertices of
the triangle $T$. Without loss of generality, we may assume that
$|AB|={\rm diam}\,T$, and, consequently, $C\ge A$ and $C\ge B$. Let
$T_A$, $T_B$, and $T_C$ be triangles obtained by joining the point
$M$ with the vertices $A$, $B$ and $C$ (see
Figure~\ref{fig:triangle_1}).

Set
$$
  m(A):=\int\!\!\!\int\limits_{T_A}(R^2-x^2-y^2)^p\,dx\,dy,
$$
$$
  m(B):=\int\!\!\!\int\limits_{T_B}(R^2-x^2-y^2)^p\,dx\,dy,
$$
$$
  m(C):=\int\!\!\!\int\limits_{T_C}(R^2-x^2-y^2)^p\,dx\,dy.
$$
Note that in the case of non-obtuse triangles
$$
  d^p(\overline{Q},T,L_p)=m(A)+m(B)+m(C),
$$
and for obtuse triangles ($C>\frac{\pi}2$) we have
$$
  d^p(\overline{Q},T,L_p)=m(A)+m(B)-m(\pi-C).
$$
Let the triangle $T_A^1$ be homothetic to the triangle $T_A$ with
the side, corresponding to the side $OB$, equal to $1$. From the
definition of $m(A)$ it follows that
\begin{equation}\label{exp_m_A}
    m(A) = \displaystyle R^{2p+2}\int\!\!\!\int\limits_{T_A^1}(1-x^2-y^2)^p\,dx\,dy  = \displaystyle 2R^{2p+2}\left(\frac A{2p+2}-l(A)\right),
\end{equation}
where
$$
  l(A)=\int\limits_{\cos{A}}^1x(1-x^2)^p\arccos{\left(\frac{\cos{A}}{x}\right)}\,dx.
$$
Note that $l(0)=0$, $l(\pi/2)=\displaystyle\frac{\pi}{4p+4}$ and
$$
  l'(A)=\gamma(p)(\sin{A})^{2p+2},
$$
where $\gamma(p)=\frac 12 B\left(p+1,\,\frac 12\right)$. 

\begin{figure}[here] 
   \centering
   \includegraphics[width=2.8in]{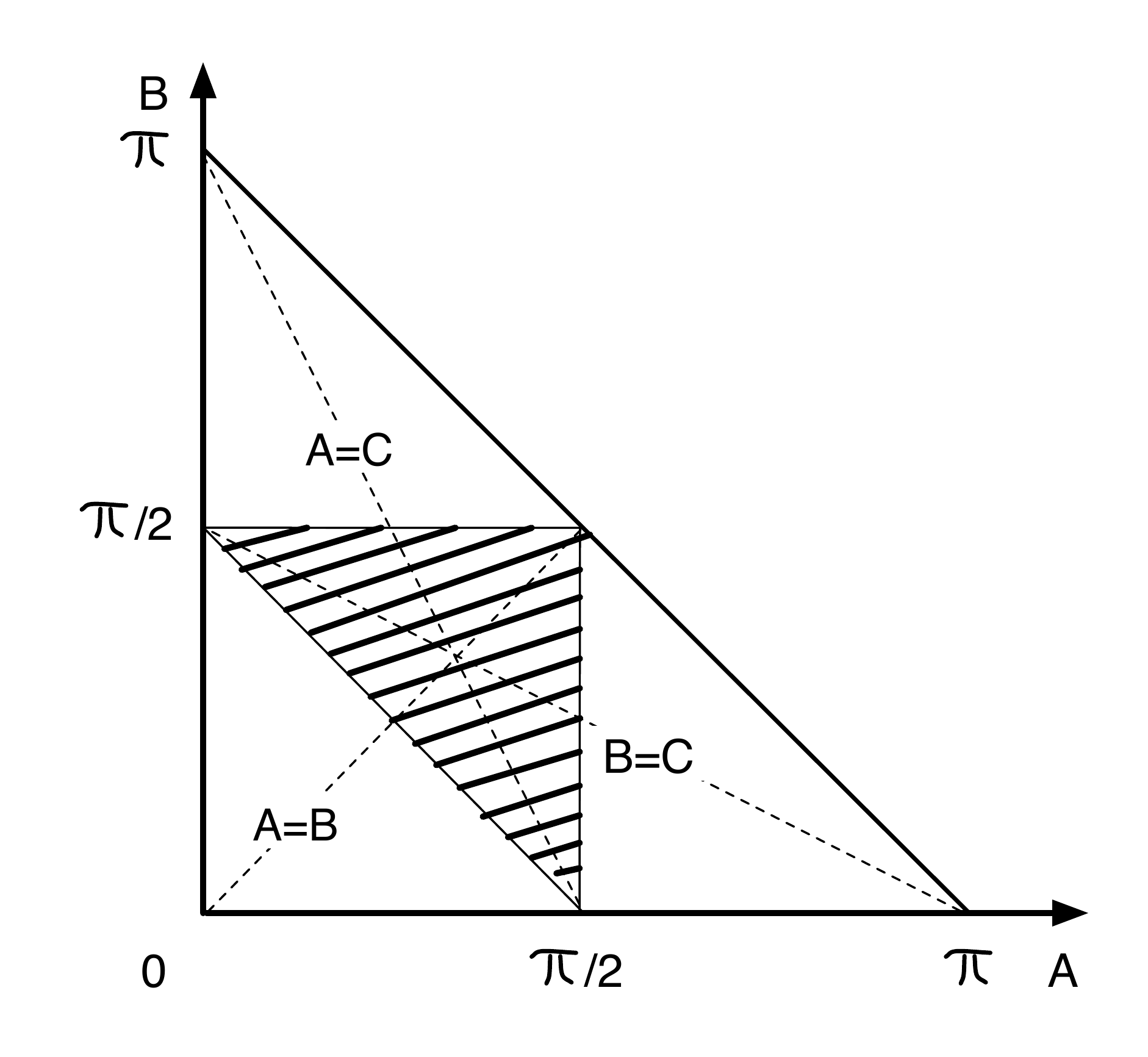} \qquad
   \includegraphics[width=3in]{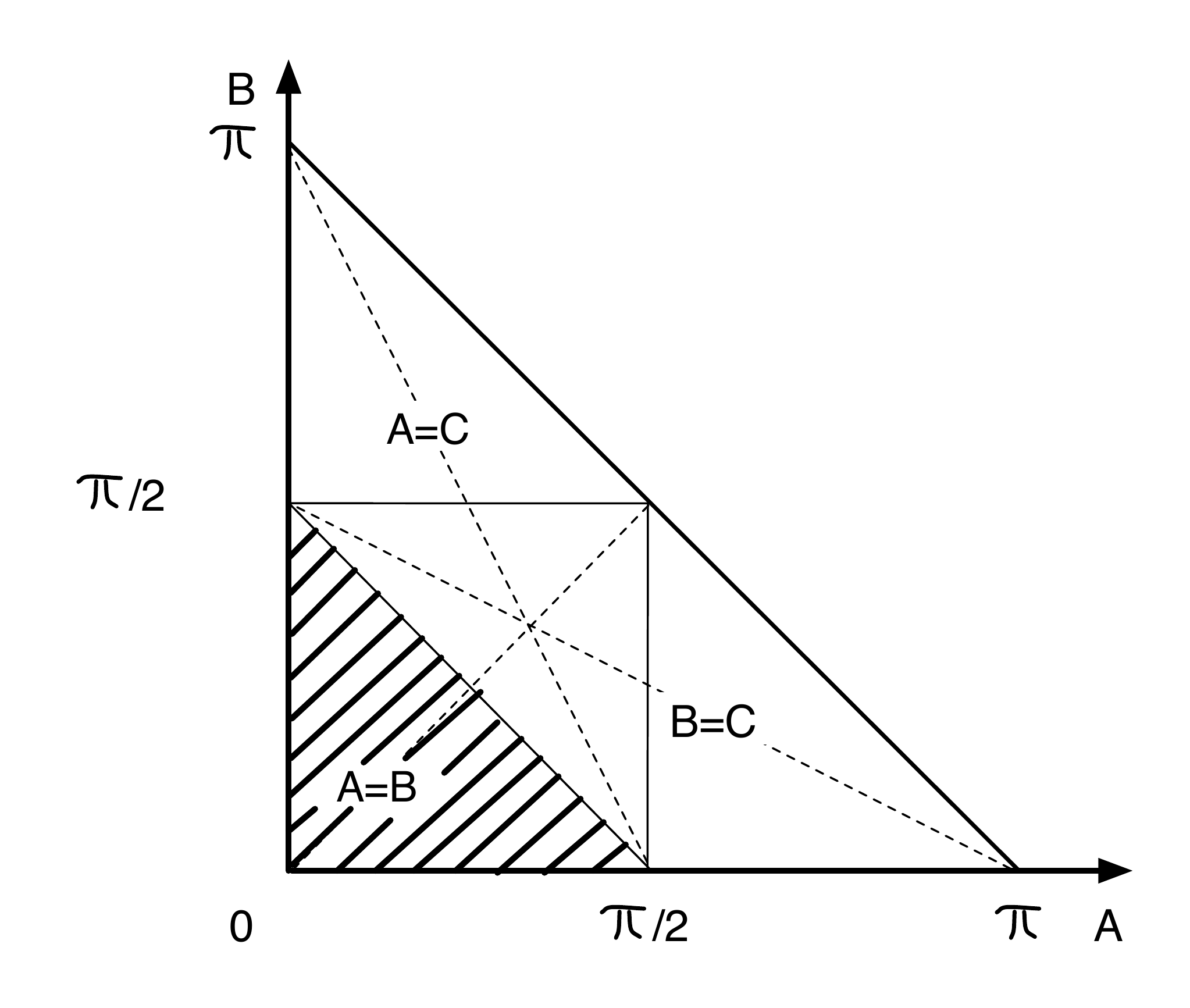}
   \caption{\it The domains $F_1$ (left) and $F_2$ (right)}
   \label{fig:reg1}
\end{figure}

Every triangle with unit area is uniquely determined by the
pair of its angles $(A,B)$, which can be considered as the point in the plane with coordinates $(A,B)$.
Let $F_1$ and $F_2$ be the sets of
points, which corresponds to the sets of all acute and all obtuse
triangles respectively (see~Figure~\ref{fig:reg1}). Obviously, $F_1$
and $F_2$ are open sets in the usual topology of $\RR^2$.

In view of Lemma~\ref{frf}, the infimum in~(\ref{inf})
is achieved at some triangle $T_1$. Let $A_1$, $B_1$ and $C_1$
be the angles of the triangle $T_1$.

First, assume that the triangle $T_1$ is acute. Taking into account
the expression~(\ref{exp_m_A}), in this case we have
$$
  d^p(\overline{Q},T,L_p)=m(A)+m(B)+m(C)=2R^{2p+2}\left[\frac{\pi}{2p+2}-l(A)-l(B)-l(C)\right].
$$
Since $|T|=1$, we obtain
$$
  R^2=\frac 2{\sin{2A}+\sin{2B}+\sin{2C}}.
$$
Hence,
\begin{equation}
  \displaystyle Z(A,B):=\frac{d^p(\overline{Q},T,L_p)}{2^{p+2}}=\displaystyle\frac{\frac{\pi}{2p+2}-l(A)-l(B)-l(C)}{(\sin{2A}+\sin{2B}+\sin{2C})^{p+1}}. \label{min_T_A_B}
\end{equation}
Therefore, by the choice of the triangle $T_1$,
$$
  Z(A_1,B_1)=\inf_{(A,B)\in F_1} Z(A,B).
$$
Since $T_1$ is an acute triangle, the necessary conditions of
extremum have to be satisfied:
\begin{equation}
  \begin{array}{rclcl}
    \displaystyle\frac{\partial Z}{\partial A} & = & \left\{\displaystyle\frac{l'(C)-l'(A)}{(\sin{2A}+\sin{2B}+\sin{2C})^{p+1}}\right. & & \\ [12pt]
                                  &   & \;\; - \;\; \left.\left.\displaystyle\frac{(2p+2)(\cos{2A}-\cos{2C})\left[\frac{\pi}{2p+2}-l(A)-l(B)-l(C)\right]}{(\sin{2A}+\sin{2B}+\sin{2C})^{p+2}}\right\}\right|_{{A=A_1\atop B=B_1}} & = & 0
  \end{array}\label{dT_dA}
\end{equation}
and
\begin{equation}
  \begin{array}{rclcl}
    \displaystyle\frac{\partial Z}{\partial B} & = & \left\{\displaystyle\frac{l'(C)-l'(B)}{(\sin{2A}+\sin{2B}+\sin{2C})^{p+1}}\right. & & \\ [12pt]
                                  &   & \;\; - \;\; \left.\left.\displaystyle\frac{(2p+2)(\cos{2B}-\cos{2C})\left[\frac{\pi}{2p+2}-l(A)-l(B)-l(C)\right]}{(\sin{2A}+\sin{2B}+\sin{2C})^{p+2}}\right\}\right|_{{A=A_1\atop B=B_1}} & = & 0
  \end{array} \label{dT_dB}
\end{equation}
From equations~(\ref{dT_dA})~--~(\ref{dT_dB}) we derive that either the angle $C_1$ is equal to one of
the other angles of the triangle $T_1$ ($C_1=A_1$ or $C_1=B_1$) or
$$
  \frac{l'(C_1)-l'(A_1)}{(\sin{C_1})^2-(\sin{A_1})^2}=\frac{l'(C_1)-l'(B_1)}{(\sin{C_1})^2-(\sin{B_1})^2}.
$$
From the last equation with the help of the expression for $l'(A)$ we have:
$$
  \frac{(\sin{C_1})^{2p+2}-(\sin{A_1})^{2p+2}}{(\sin{C_1})^2-(\sin{A_1})^2}=\frac{(\sin{C_1})^{2p+2}-(\sin{B_1})^{2p+2}}{(\sin{C_1})^2-(\sin{B_1})^2}
$$
which can be true only if $A_1=B_1$.
Therefore, if $T_1$ is acute then it is isosceles.

Now, if we assume that $T_1$ is obtuse then we similarly obtain that
$T_1$ should be isosceles. Therefore, the triangle $T_1$
solving~(\ref{inf}) must be either right or isosceles.

Now let us consider the case of right triangles. Then $C=\frac{\pi}2$. We shall show that the
triangle solving the problem
$$
  d(\overline{Q},T,L_p)\to\inf, \quad |T|=1,\; C=\frac{\pi}2,
$$
is isosceles. Indeed, in this case $B=\frac{\pi}2 - A$ and
$$
  d^p(\overline{Q},T,L_p)=m(A)+m\left(\frac{\pi}2 - A\right) + m\left(\frac{\pi}2\right) = 2\cdot\displaystyle\frac{\displaystyle\frac{\pi}{4p+4}-l(A)-l\left(\frac{\pi}2-A\right)}{(\sin{2A})^{p+1}}.
$$

\begin{lemma}
The function
\begin{equation}\label{LA}
  L(A):=\displaystyle\frac{\displaystyle\frac{\pi}{4p+4}-l(A)-l\left(\frac{\pi}2-A\right)}{(\sin{2A})^{p+1}}
\end{equation}
is non-increasing on the interval $(0,\pi/4]$.
\end{lemma}

{\bf {Proof.}}
Let us consider the derivative of the function $L$:
$$
  L'(A)=\frac{[l'(\pi/2-A)-l'(A)]\sin{2A}-(2p+2)\cos{2A}\left[\displaystyle\frac{\pi}{4p+4}-l(A)-l(\pi/2-A)\right]}{(\sin{2A})^{p+2}}.
$$
Taking into account the expression for $l'(A)$, in order to show
that $L'(A)\le0$ on the interval $(0,\pi/4)$ it suffices to prove
that
$$
  r(A):=\gamma(p)[(\cos{A})^{2p+2}-(\sin{A})^{2p+2}]\tan{2A}-(2p+2)\left[\frac{\pi}{4p+4}-l(A)-l(\pi/2-A)\right]\le0
$$
on the interval $A\in(0,\pi/4)$.

Let us consider the derivative $r'(A)$:
$$
  \begin{array}{rcl}
     r'(A) & = & \gamma(p)\Big\{\displaystyle\frac 2{(\cos{2A})^2}[(\cos{A})^{2p+2}-(\sin{A})^{2p+2}] \\ [12pt]
           &   & \begin{array}{ll}
                 \quad\quad\; - & (2p+2){\tan}\,2A\sin{A}\cos{A}[(\sin{A})^{2p}+(\cos{A})^{2p}] \\ [7pt]
                 \quad\quad\; - & (2p+2)[(\cos{A})^{2p+2}-(\sin{A})^{2p+2}]\Big\},
                 \end{array}
  \end{array}
$$
where $\gamma(p)=\frac 12 B\left(p+1,\,\frac 12\right)$.

Set $t={\rm tan}\,A$, and then $t\in(0,1]$. Obviously,
$$
  r'(A)=\frac{\gamma(p)z(t)}{(\cos{A})^{2p+2}},
$$
where
$$
    z(t) 
             =pt^{2p+6}-2t^{2p+4}-(p+2)t^{2p+2}+(p+2)t^4+2t^2-p.
$$

To prove that $r'(A)\le 0$ for all $A\in\left(0,\frac{\pi}{4}\right]$
we shall need the following proposition.

{\bf Proposition 1. \it Let $p\in(0,1)$. Then the function
$$
  z(t)=pt^{2p+6}-2t^{2p+4}-(p+2)t^{2p+2}+(p+2)t^4+2t^2-p
$$
is non-positive for all $t\in[0,1]$. }

Thus, we obtain that $r'(A)\le0$ for all
$A\in\displaystyle\left(0,\frac{\pi}{4}\right]$. Consequently, the
function $r(A)$ is non-increasing, and hence $r(A)\le r(+0)=0$. Hence,
we have proved that the function $L(A)$ is non-increasing on the
interval $(0, \pi/4]$, i.e.
$$
  L(A)\ge L(\pi/4)
$$
for all $A\in (0, \pi/4)$. $\square$

Note that it can be shown similarly that the function $L(A)$ (see~(\ref{LA}))
is non-decreasing on the interval $\displaystyle\left[\frac{\pi}4,\frac{\pi}2\right)$.
Therefore,
$$
  \inf\limits_{|T|=1,\,C=\frac{\pi}2}d(\overline{Q},T,L_p) = d(\overline{Q}, T_2,L_p),
$$
where $T_2$ is the isosceles right triangle.

Hence, we have proved that the triangle $T_1$ must be isosceles.
Without loss of generality we may assume that $A=B$. Let us show
that $T_1$ is equilateral.
To this end let us consider the following problem
\begin{equation}\label{A=B}
  d(\overline{Q},T,L_p)\to\inf,\quad |T|=1,\; A=B.
\end{equation}

Let us consider now the case when the triangle $T$ is isosceles non-obtuse triangle. Then
$A\in\displaystyle\left[\frac{\pi}{4},\frac{\pi}{2}\right)$, and
$$
  d^p(\overline{Q},T,L_p) = 2m(A) + m(\pi-2A) = 2^{p+2}\cdot\displaystyle\frac{\frac{\pi}{2p+2}-2l(A)-l(\pi-2A)}{(2\sin{2A}-\sin{4A})^{p+1}}.
$$

\begin{lemma}
The function
$$
  S(A)= \displaystyle\frac{\frac{\pi}{2p+2}-2l(A)-l(\pi-2A)}{(2\sin{2A}-\sin{4A})^{p+1}}
$$
is non-increasing for $A\in[\pi/4,\pi/3]$ and is non-decreasing
for $A\in[\pi/3,\pi/2)$.
\end{lemma}

{\bf {Proof.}}
Let us consider the derivative of the function $S$
\begin{equation}
  \begin{array}{rcl}
    S'(A) & = & \displaystyle\frac{8\cos{A}(\sin{A})^3[2l'(\pi-2A)-2l'(A)]}{(8\cos{A}(\sin{A})^3)^{p+2}} \\[12pt]
          &   & \;\; - \;\; \displaystyle\frac{8(p+1)\sin{A}\sin{3A}\left[\frac{\pi}{2p+2}-2l(A)-l(\pi-2A)\right]}{(8\cos{A}(\sin{A})^3)^{p+2}}.
  \end{array} \label{dS_dA}
\end{equation}
Obviously, the denominator of the right hand side of~(\ref{dS_dA})
is positive in the considered region. In what follows we shall
consider only the numerator of ~(\ref{dS_dA}). Then showing that
$S'(A)\le 0$, when $A\in[\pi/4,\pi/3]$ and $S'(A)\ge 0$ when
$A\in[\pi/3,\pi/2)$, is equivalent to proving the inequality:
$$
  \begin{array}{rcl}
    q(A) & := & 2\gamma(p)\displaystyle\frac{\cos{A}\sin{A}}{3-4(\sin{A})^2}\Big[(\sin{2A})^{2p+2}-(\sin{A})^{2p+2}\Big] \\ [12pt]
         &    & \;\; - \;\; (p+1)\left[\displaystyle\frac{\pi}{2p+2}-2l(A)-l(\pi-2A)\right]\le 0
  \end{array}
$$
for all $A\in[\pi/4, \pi/2)$.

It is easy to see that $q(\pi/2-0)=0$ and
$$
  \begin{array}{rcl}
    q'(A) &   =   & \displaystyle\frac{2\gamma(p)}{\left[3-4(\sin{A})^2\right]^2}\Big\{\left[(\sin{2A})^{2p+2}-(\sin{A})^{2p+2}\right] \\ [7pt]
          &       &  \begin{array}{ll}
                      \qquad\qquad\qquad\qquad\times & \Big[\cos{2A}\Big(3-4(\sin{A})^2\Big)+8(\cos{A})^2(\sin{A})^2 \\ [7pt]
                                                     & \;\; - \;\; (p+1)\Big(3-4(\sin{A})^2\Big)^2\Big] \\ [7pt]
                      \qquad\qquad\qquad\qquad   +   & 2(p+1)\cos{A}\sin{A}\Big[3-4(\sin{A})^2\Big] \\[7pt]
                                                     & \;\;\times \;\; \left[2\cos{2A}(\sin{2A})^{2p+1}-\cos{A}(\sin{A})^{2p+1}\right]\Big\}.
                     \end{array}
  \end{array}
$$
Obviously,
$$
  q'(A) = \frac{\gamma(p)}{\left[3-4(\sin{A})^2\right]^2}z(t), \quad t=2\cos{A}\in(0,\sqrt{2}],
$$
where
$$
  \begin{array}{rcl}
    z(t) & := & (t^{2p+2}-1)\cdot\Big[(t^2-2)(t^2-1)+t^2(4-t^2)-2(p+1)(t^2-1)^2\Big] \\ [7pt]
         &    & \qquad\qquad\qquad\;\; + \;\; (p+1)t(t^2-1)\Big[2(t^2-2)t^{2p+1}-t\Big] \\ [7pt]
         &  = & -(2p+1)t^{2p+4}+2(p+2)t^{2p+2}+(p+1)t^4-(3p+4)t^2+2p.
  \end{array}
$$

We shall need the following proposition to show that $z(t)$ changes its sign exactly once.

{\bf Proposition 2. \it
 Let $p\in(0,1)$. Then the function
$$
  z(t)=-(2p+1)t^{2p+4}+2(p+2)t^{2p+2}+(p+1)t^4-(3p+4)t^2+2p
$$
has exactly one point of sign change (from positive to negative) on the segment $[0,2]$,
and this point is located inside the interval $(0,1)$.
}

Therefore, $z(t)$ changes its sign on the segment $[0,\sqrt{2}]$
exactly once from positive to negative. Since variable $A$ increases
when variable $t$ decreases, the function $q'(A)$ changes sign
exactly once as well, from negative to positive. Hence,
$$
  q(A)\le\min\left\{q(\pi/4);\,q(\pi/2-0)\right\}=0,
$$
whenever
$q(\pi/4)=\gamma(p)(1-2^{-p-1})-(p+1)\left[\frac{\pi}{4p+4}-2l(\pi/4)\right]\le
0$ (this will be shown in the proof of Lemma 13). Hence, we have
proved that $S'(A)\le 0$ for $A\in[\pi/4, \pi/3]$, and $S'(A)\ge 0$
for $A\in[\pi/3, \pi/2)$. $\square$

Let us turn now to the case of isosceles obtuse triangles $T$.
In this case $A\in\left(0,\frac{\pi}{4}\right)$, and
$$
  d^p(\overline{Q},T,L_p) = 2m(A) - m(2A) = 2^{p+2}\cdot\displaystyle\frac{-2l(A)+l(2A)}{(2\sin{2A}-\sin{4A})^{p+1}}.
$$

\begin{lemma}
The function
$$
  {\widetilde S}(A) = \displaystyle\frac{-2l(A)+l(2A)}{(2\sin{2A}-\sin{4A})^{p+1}}
$$
is non-increasing on the interval $\left(0,\frac{\pi}{4}\right]$.

\end{lemma}
{\bf Proof.} Let us
 consider the derivative of the function $\widetilde{S}$
$$
  \begin{array}{rcl}
    {\widetilde S}'(A) & = & \displaystyle\frac{8\cos{A}(\sin{A})^3[2l'(2A)-2l'(A)]}{(8\cos{A}(\sin{A})^3)^{p+2}} \\[12pt]
                         & & \;\;-\;\;\displaystyle\frac{8(p+1)\sin{A}\sin{3A}[-2l(A)+l(2A)]}{(8\cos{A}(\sin{A})^3)^{p+2}}.
  \end{array}
$$
We shall show that ${\widetilde S}'(A)\le 0$ when $A\in(0,\,\pi/4]$
or, equivalently, that ${\tilde q}(A)\le0$ for all $A\in(0,\,\pi/4]$
where
$$
  \begin{array}{rcl}
    {\tilde q}(A) & := & 2\gamma(p)\displaystyle\frac{\cos{A}\sin{A}}{3-4(\sin{A})^2}\Big[(\sin{2A})^{2p+2}-(\sin{A})^{2p+2}\Big] \\ [12pt]
         &   & \;\;-\;\; (p+1)[-2l(A)+l(2A)].
  \end{array}
$$

The derivative of ${\tilde q}(A)$ can be written as:
$$
  \begin{array}{rcl}
    {\tilde q}'(A) &   =   & \displaystyle\frac{2\gamma(p)}{\left[3-4(\sin{A})^2\right]^2}\Big\{\left[(\sin{2A})^{2p+2}-(\sin{A})^{2p+2}\right] \\ [7pt]
          &&    \begin{array}{ll}
                  \qquad\qquad\qquad\quad \times & \Big[\cos{2A}\Big(3-4(\sin{A})^2\Big)+8(\cos{A})^2(\sin{A})^2\\[7pt]
                                                  & \;\;- \;\;(p+1)\Big(3-4(\sin{A})^2\Big)^2\Big] \\ [7pt]
                  \qquad\qquad\qquad\quad    +   & 2(p+1)\cos{A}\sin{A}\Big[3-4(\sin{A})^2\Big] \\ [7pt]
                                                  & \;\;\times\;\;\left[2\cos{2A}(\sin{2A})^{2p+1}-\cos{A}(\sin{A})^{2p+1}\right]\Big\}
                \end{array} \\ [7pt]
          &   =   & \displaystyle\frac{\gamma(p)}{\left[3-4(\sin{A})^2\right]^2}z(t),
  \end{array}
$$
where $t=2\cos{A}\in[\sqrt{2},\,2)$ and
$$
  z(t) = -(2p+1)t^{2p+4}+2(p+2)t^{2p+2}+(p+1)t^4-(3p+4)t^2+2p.
$$
In view of Proposition~2, it follows that $z(t)\le 0$ for all
$t\in[\sqrt{2},\,2)$. Hence, ${\tilde q}'(A)\le 0$ on the segment
$(0,\,\pi/4]$. Therefore,  ${\tilde q}(A)\le {\tilde q}(+0)=0$ (in
particular, $q(\pi/4)={\tilde q}(\pi/4)\le 0$) and ${\tilde S}(A)$
is non-increasing on the segment $(0,\,\pi/4]$. $\square$

Combining Lemmas~12 and~13 we conclude that
$$
  \inf_{|T|=1,\,A=B}d(\overline{Q},T,L_p)=d(\overline{Q},T_3,L_p),
$$
where $T_3$ is equilateral triangle. Moreover, $T_3$ is the only
triangle which solves~(\ref{A=B}). Since the triangle $T_1$
solves~(\ref{inf}) and is isosceles, we obtain that $T_1=T_3$,
which finishes the proof.

Department of Mathematical Analysis \\
Dnepropetrovsk National University \\
Institute of Applied Mathematics and Mechanics of NAS of Ukraine \\
pr. Gagarina, 72, \\
Dnepropetrovsk, UKRAINE, 49050 \\
Phone: +380567609461\\
Email: babenko.vladislav@gmail.com \\

Department of Mathematics and Statistics\\
Sam Houston State University\\
Box 2206\\
Huntsville, TX, USA 77340-2206\\
Phone: 936.294.4884\\
Fax: 936.294.1882\\
Email: babenko@shsu.edu\\

Department of Mathematical Analysis \\
Dnepropetrovsk National University \\
pr. Gagarina, 72 \\
Dnepropetrovsk, UKRAINE, 49050 \\
Email: dmitriy.skorokhodov@gmail.com\\

\end{document}